\documentclass[a4paper,12pt]{amsart}

\usepackage{amsfonts}
\usepackage{amsmath}
\usepackage{amssymb}
\usepackage{mathrsfs}
\usepackage{hyperref}
\usepackage{graphicx}

\usepackage[usenames]{color}

\numberwithin{equation}{section}
\theoremstyle{plain}
\newtheorem{thm}{Theorem}[section]

\theoremstyle{definition}
\newtheorem{defi}[thm]{Definition}

\newtheorem{prm}[thm]{Problem}
\def\L{{\mathcal L}}
\def\p#1{{\left({#1}\right)}}

\def\R{\mathcal R}
\def\G{{\mathbb G}}
\def\H{\mathcal H}
\newcommand{\Dcal}{\mathcal D}

\begin{document}

\title[Inverse source problems for positive operators]
{Inverse source problems for positive operators. I. Hypoelliptic diffusion and subdiffusion equations}

\author[Michael Ruzhansky]{Michael Ruzhansky}
\address{
  Michael Ruzhansky:
  \endgraf
Department of Mathematics: Analysis,
Logic and Discrete Mathematics
  \endgraf
Ghent University, Belgium
  \endgraf
 and
  \endgraf
 School of Mathematical Sciences
 \endgraf
Queen Mary University of London
\endgraf
United Kingdom
\endgraf
  {\it E-mail address} {\rm michael.ruzhansky@ugent.be}
 }

\author[Niyaz Tokmagambetov]{Niyaz Tokmagambetov}
\address{
  Niyaz Tokmagambetov:
  \endgraf
    Al--Farabi Kazakh National University
  \endgraf
  71 Al--Farabi ave., Almaty, Kazakhstan
  \endgraf
  and
  \endgraf
Department of Mathematics: Analysis,
Logic and Discrete Mathematics
  \endgraf
Ghent University, Belgium
  \endgraf
and
  \endgraf
  Institute of Mathematics and Mathematical Modeling
  \endgraf
  125 Pushkin str., Almaty, Kazakhstan
  \endgraf
  {\it E-mail address} {\rm niyaz.tokmagambetov@ugent.be}
 }

\author[Berikbol T. Torebek]{Berikbol T. Torebek}
\address{
  Berikbol T. Torebek:
  \endgraf
    Al--Farabi Kazakh National University
  \endgraf
  71 Al--Farabi ave., Almaty, Kazakhstan
  \endgraf
  and
  \endgraf
Department of Mathematics: Analysis,
Logic and Discrete Mathematics
  \endgraf
Ghent University, Belgium
  \endgraf
and
  \endgraf
  Institute of Mathematics and Mathematical Modeling
  \endgraf
  125 Pushkin str., Almaty, Kazakhstan
  \endgraf
  {\it E-mail address} {\rm berikbol.torebek@ugent.be}
  }

\thanks{The first author was supported in parts by the FWO Odysseus Project, EPSRC grant EP/R003025/1 and by the Leverhulme Grant RPG-2017-151. The second author was supported by the Ministry of Education and Science of the Republic of Kazakhstan Grant AP05130994. The third author was supported by the Ministry of Education and Science of the Republic of Kazakhstan Grant AP05131756. No new data was collected or generated during the course of research.}

\date{\today}


\begin{abstract}
A class of inverse problems for restoring the right-hand side of a parabolic equation for a large class of positive operators with discrete spectrum is considered. The results on existence and uniqueness of solutions of these problems as well as on the fractional time diffusion (subdiffusion) equations are presented. Consequently, the obtained results are applied for the similar inverse problems for a large class of subelliptic diffusion and subdiffusion equations (with continuous spectrum). Such problems are modelled by using general homogeneous left-invariant hypoelliptic operators on general graded Lie groups. A list of examples is discussed, including Sturm-Liouville problems, differential models with involution, fractional Sturm-Liouville operators, harmonic and anharmonic oscillators, Landau Hamiltonians, fractional Laplacians, and harmonic and anharmonic operators on the Heisenberg group. The rod cooling problem for the diffusion with involution is modelled numerically, showing how to find a ``cooling function'', and how the involution normally slows down the cooling speed of the rod.
\end{abstract}

\maketitle
\tableofcontents
\section{Introduction}
Many instances are known in which the practical needs lead to the problems of determining the coefficients or the right-hand-side of a differential equation from some available data about the solution.  These are called the inverse problems of mathematical  physics.   Inverse  problems  arise  in  various  areas  of  human  activity  such  as  seismology,  mineral  exploration,  biology,  medicine,  quality  control  of industrial goods,  etc.  All these circumstances place inverse problems among the important problems of modern mathematics.

The first purpose of this paper is to study inverse problems for the heat equation. We consider the heat equation
\begin{equation}\label{1}
u_t(x,t)+\mathcal{L}u(x,t)=f(x),
\end{equation}
with Cauchy condition
\begin{equation}
\label{2}
u(x,0)=\varphi(x), x\in \bar{Q},
\end{equation}
for $(x,t)\in\Omega_t=\{x\in Q\subset \mathbb{R}^d,\,d\geq 1,\,\, t\in[0,T]\}$, where $\mathcal L$ is a linear self-adjoint positive operator with a discrete spectrum $\{\lambda_\xi>0:\xi\in \mathcal I\}$ on a separable Hilbert space $\mathcal H$. Here $Q$ is a bounded domain with a smooth boundary or unbounded domain. Respectively $\lambda_\xi$, the operator $\mathcal L$ has the system of orthonormal eigenfunctions $\{e_\xi:\xi\in \mathcal I\}$ on the separable Hilbert function space $\mathcal H$.

The problem of determination of temperature at interior points of a region when the initial and boundary conditions along with diffusion source term are specified are known as direct diffusion conduction problems. In many physical problems, determination of coefficients or right hand side (the source term, in case of the diffusion equation) in a differential equation from some available information is
required; these problems are known as inverse problems. These kind of problems are ill posed in the sense of Hadamard. A number of articles address the analytical and numerical solvability of the inverse problems for the diffusion and anomalous diffusion equations (see \cite{CD98, CNYY09, JR15, Kali1, Oraz, Oraz1, torebek, Fur, Ismailov, Kir, Kir1, Nguyen, ZX11, ZW13, WYH13} and references therein).

The setting of a general operator $\mathcal{L}$ as in this paper allows one to include many models. A number of physical examples are discussed in Section \ref{SEC:examples}, including Sturm-Liouville problems, differential models with involution, fractional Sturm-Liouville operators, harmonic and anharmonic oscillators, Landau Hamiltonians, fractional Laplacians, and harmonic and anharmonic operators on the Heisenberg group.

Section \ref{InvFr} is dedicated to finding the couple of functions $\left(u( {x,t}
), f( x )\right)$ satisfying the equation
\begin{equation}\label{I-R: 2.1a}
\mathcal{D}_{0+,t}^{\alpha} u \left( {x,t} \right) + \mathcal{L} u(x, t) = f\left( x \right),
\end{equation}
in the domain $(x,t)\in\Omega_t=\{x\in Q\subset \mathbb{R}^d,\,d\geq 1,\,\, t\in[0,T]\},$
under the conditions \begin{equation*}u\left(
{x,0} \right) = \varphi \left( x \right),\  x \in \bar Q,\end{equation*} \begin{equation*}u\left( {x,T} \right) =
\psi \left( x \right),\ x \in \bar Q,
\end{equation*} where $\mathcal{D}_{0+,t}^{\alpha}$ is a Caputo fractional derivative of order $0<\alpha\leq 1,$ $\varphi \left( x \right)$ and $\psi
\left( x \right)$ are sufficiently smooth functions, and $\mathcal L$ is a linear positive operator with a discrete spectrum.

In many contexts, for example, in the sub-Riemannian settings, the (fractional) time diffusion equation for subelliptic operators arises naturally. However, such operators have a non-discrete, but continuous spectrum. Following Rothschild and Stein \cite{RS76} and further developments, many of such operators can be modelled by the so-called Rockland operators (homogeneous left-invariant hypoelliptic differential operators) on graded Lie groups. {\em By using the proceeding analysis, we will employ the group Fourier transform to reduce such problems to those with the discrete spectrum, for which the above established results would be applicable.}

Thus, let $\R$ be
a positive self--adjoint Rockland operator acting on $L^{2}(\G)$, where $\G$ is a graded Lie group of homogeneous dimension $Q\geq 3$. In a domain $\Omega=\{(t,x):\,(0,T)\times \G\}$ we seek a couple of functions $\left( u(t), f \right)$ satisfying the equation
\begin{equation}\label{I-R: 2.1}
\mathcal{D}_{0+,t}^{\alpha} u \left( {x,t} \right) + \R u(x, t) = f\left( x \right),
\end{equation}
under the conditions
\begin{equation}\label{I-R: 2.2}
u\left({x,0} \right) = \varphi \left( x \right),\  x \in \G,
\end{equation}
\begin{equation}\label{I-R: 2.2*}
u\left( {x,T} \right) =
\psi \left( x \right),\ x \in \G,
\end{equation}
where $\mathcal{D}_{0+,t}^{\alpha}$ is a Caputo fractional derivative, $0 < \alpha \leq 1$.

We seek a solution $u\in C([0,T]; L^{2}(\G))$ of the problem \eqref{I-R: 2.1}--\eqref{I-R: 2.2*} such that $\mathcal{D}_{0+,t}^{\alpha} u\in C([0,T]; L^{2}(\G))$, $\R u\in C([0,T]; L^{2}(\G))$, and $f\in C([0,T]; L^{2}(\G))$. We are able to solve this problem by employing the natural global Fourier analysis on $\G$, allowing one to use the solution to the problem \eqref{I-R: 2.1a} applied to the infinitesimal representations of the operator $\R$, which is in turn known to have the discrete spectrum.

We note that the inverse source problems for Rockland operators \eqref{I-R: 2.1}--\eqref{I-R: 2.2*} cover also the Heisenberg case, namely,
the equation
\begin{equation*}
\mathcal{D}_{0+,t}^{\alpha} u \left( {x,t} \right) - \Delta_{\mathbb{H}^n} u(x, t) = f\left( x \right),
\end{equation*}
with the conditions
\begin{equation*}
u\left({x,0} \right) = \varphi \left( x \right),\  x \in \mathbb{H}^n,
\end{equation*}
\begin{equation*}
u\left( {x,T} \right) =
\psi \left( x \right),\ x \in \mathbb{H}^n,
\end{equation*}
where $\mathcal{D}_{0+,t}^{\alpha}$ is a Caputo fractional derivative, $0 < \alpha \leq 1$. Here, $\Delta_{\mathbb{H}^n}$ is the sub-Laplacian on the Heisenberg group $\mathbb{H}^n$. As the physical application, the relevance of the Heisenberg group $\mathbb{H}^n$ for quantum mechanics has been established. In 1931 Weyl \cite{W31} recognized that the Heisenberg algebra generated by the momentum and position operators originates from the representation of the Lie algebra associated with the corresponding group, namely, the Heisenberg group, or the Weyl group as the physicists call it. For the recent studies of the heat equation on the Heisenberg groups, we refer to the publications \cite{BC17, DM05, Jui14, TW18, RTT19}.

Thus, let us briefly summarise the results of this paper:
\begin{itemize}
\item Existence and uniqueness for the inverse diffusion problem
$$u_t(x,t)+\mathcal{L}u(x,t)=f(x),$$
for general positive operators $\mathcal{L}$ with discrete spectrum and the basis of eigenfunctions.
\item Existence and uniqueness for the inverse time-fractional subdiffusion problem
$$\mathcal{D}_{0+,t}^{\alpha} u \left( {x,t} \right) + \mathcal{L} u(x, t) = f( x ).$$
\item Existence and uniqueness for the inverse diffusion and subdiffusion problems
$$\mathcal{D}_{0+,t}^{\alpha} u \left( {x,t} \right) + \R u(x, t) = f( x ),$$
for $0 < \alpha \leq 1$, and for general homogeneous left-invariant hypoelliptic differential operators $\R$ on graded Lie groups.
\item Numerical analysis of the obtained formulae in the case of the inverse heat problems for differential operators with involution.
\end{itemize}

\section{Inverse problem for the heat equation}

\subsection{Statement of the problem}
The Section is concerned with inverse problem for the heat equation \eqref{1}. We obtain existence and uniqueness results for this problem, based on the $\L$--Fourier method. An introduction and some basic definitions of the $\L$--Fourier analysis are given in \cite{KRT17, DRT17, RT16,  RT17a, RT18b}.

\begin{prm}
\label{PR-01}
Find the couple of functions $(u(x,t),f(x))$ satisfying the Cauchy problem \eqref{1}-\eqref{2}, under the condition
\begin{equation}
\label{3}
u(x,T)=\psi(x), \,\,\, x\in \bar{Q}.
\end{equation}
\end{prm}

A generalised solution of Problem \ref{PR-01} is the pair of functions $(u(x,t),f(x)),$ where $u=u(x,t)\in C^1([0,T],\mathcal H)\cap C([0,T],\mathcal H^{1})$, and $f=f(x)\in \mathcal H$. Here $\mathcal H^{s}$ is the closure of $\mathcal{H}_\mathcal{L}^\infty$ (see, Appendix \ref{App 01}) under the norm
$$
\|v\|_{\mathcal H^{s}}:=\|(I+\L)^{s}v\|_{\mathcal H},
$$
for all $s\in\mathbb R_{+}$.

For the considered Problem \ref{PR-01}, the following theorem holds true.

\begin{thm}
\label{Th-01}
Let $\varphi, \psi\in \mathcal H^{1}$. Then the generalised solution of Problem \ref{PR-01} exists, is unique, and can be written in the form
$$u(x,t)=\varphi(x)+\sum_{\xi\in \mathcal I}\frac{[(\varphi,e_\xi)_{\mathcal H}-(\psi,e_\xi)_{\mathcal H}](e^{-\lambda_\xi t}-1)e_\xi(x)}{(1-e^{-\lambda_\xi T}) },$$

$$
f(x)=\mathcal L\varphi(x)-\sum_{\xi\in \mathcal I}\frac{[(\mathcal L\varphi,e_\xi)_{\mathcal H}-(\mathcal L\psi,e_\xi)_{\mathcal H}]e_\xi(x)}{(1-e^{-\lambda_\xi T})}.
$$
\end{thm}

\subsection{Proof of the existence result}
We want to find a generalised solution by the Fourier method, we have the eigenvalues $\lambda_\xi$ and  eigenfunctions system $e_\xi(x)$  of operator $\mathcal L$ on the space $\mathcal H$. Eigenfunctions system $e_\xi(x)$ is an orthonormal basis in $\mathcal H$, the functions $u(x,t)$ and $f(x)$ can be expanded as follows:
\begin{equation}\label{4} u(x,t)=\sum_{\xi\in \mathcal I}u_\xi(t)e_\xi(x),\end{equation} and
\begin{equation}\label{5} f(x)=\sum_{\xi\in \mathcal I}f_\xi e_\xi(x),\end{equation}
where $f_\xi,u_\xi(t)$ are unknown. Substituting Equations \eqref{4} and \eqref{5} into Equation \eqref{1}, we obtain the following equation for the function $u_\xi(t)$ and the constant $f_\xi$:
$$u^{'}_{\xi}(t)+\lambda_\xi u_\xi(t)=f_\xi.$$
Solving this equation, we obtain
$$u_\xi(t)=\frac{f_\xi}{\lambda_\xi}+C_\xi e^{-\lambda_\xi t},$$ where the constants $C_\xi,f_\xi $ are unknown. To find these constants, we use conditions \eqref{2} and \eqref{3}. Let $\varphi_\xi, \psi_\xi$ be the coefficients of the expansions of $\varphi(x)$ and $\psi(x)$

$$\varphi_\xi=(\varphi,e_\xi)_{\mathcal H},$$ and
$$\psi_\xi=(\psi,e_\xi)_{\mathcal H}.$$
We first find $C_\xi$:$$u_\xi(0)=\frac{f_\xi}{\lambda_\xi}+C_\xi=\varphi_\xi,$$
and
$$
u_\xi(T)=\frac{f_\xi}{\lambda_\xi}+C_\xi e^{-\lambda_\xi T}=\psi_\xi,
$$
then $$C_\xi=\frac{\varphi_\xi-\psi_\xi}{1-e^{-\lambda_\xi T}}.$$
The constant $f_\xi$ is represented as $$f_\xi=\lambda_\xi \varphi_\xi-\lambda_\xi C_\xi.$$

Substituting $f_\xi,u_\xi(t)$ into formulas \eqref{4} and \eqref{5}, we find $$u(x,t)=\varphi(x)+\sum_{\xi\in \mathcal I}C_\xi(e^{-\lambda_\xi t}-1)e_\xi(x).$$

Using the self-adjointness property of the operator $\mathcal{L},$ we have $$(\mathcal{L}\varphi,e_\xi)_{\mathcal H}=(\varphi,\mathcal{L}e_\xi)_{\mathcal H}$$ we know that $\mathcal{L}e_\xi=\lambda_\xi e_\xi$ using this we obtain
$$(\varphi,e_\xi)_{\mathcal H}=\frac{(\mathcal{L}\varphi,e_\xi)_{\mathcal H}}{\lambda_\xi},$$ and for $\psi(x)$ we can write it in a similar way. Substituting these equality into formula of $C_\xi$ we can get that $$C_\xi=\frac{(\mathcal L\varphi,e_\xi)_{\mathcal H}-(\mathcal L\psi,e_\xi)_{\mathcal H}}{\lambda_\xi(1-e^{-\lambda_\xi T})}.$$
Then $$u(x,t)=\varphi(x)+\sum_{\xi\in \mathcal I}\frac{[(\mathcal L\varphi,e_\xi)_{\mathcal H}-(\mathcal L\psi,e_\xi)_{\mathcal H}](e^{-\lambda_\xi t}-1)e_\xi(x)}{\lambda_\xi(1-e^{-\lambda_\xi T}) }.$$
Similarly, $$f(x)=\mathcal L\varphi(x)-\sum_{\xi\in \mathcal I}\frac{[(\mathcal L\varphi,e_\xi)_{\mathcal H}-(\mathcal L\psi,e_\xi)_{\mathcal H}]e_\xi(x)}{(1-e^{-\lambda_\xi T})}.$$
Now, for the convergence of the series, using $\|e_\xi\|_{\mathcal H}=1$ and $\varphi, \psi\in \mathcal H^{1}$, we have the following estimate
\begin{align*}
\max_{t\in [0,T]}\|u(x,t)\|_{\mathcal H}^{2}&\leq C\|\varphi(x)\|_{\mathcal H}^{2}+C\max_{t\in [0,T]}\sum_{\xi\in \mathcal I}\frac{|(\varphi,e_\xi)_{\mathcal H}-(\psi,e_\xi)_{\mathcal H}|^{2}\cdot|e^{-\lambda_\xi t}-1|^{2}\cdot\|e_\xi\|_{\mathcal H}^{2}}{|1-e^{-\lambda_\xi T}|^{2} }\\&
\leq C\|\varphi(x)\|_{\mathcal H}^{2}+C\sum_{\xi\in \mathcal I}|(\varphi,e_\xi)_{\mathcal H}|^{2}+|(\psi,e_\xi)_{\mathcal H}|^{2}<\infty, \,\,\,C=const>0,
\end{align*} and
\begin{align*}
\max_{t\in [0,T]}\|u_t(x,t)\|_{\mathcal H}|^{2}&\leq\max_{t\in [0,T]}\sum_{\xi\in \mathcal I}\frac{|(\mathcal L\varphi,e_\xi)_{\mathcal H}-(\mathcal L\psi,e_\xi)_{\mathcal H}|^{2}\cdot|-\lambda_\xi\cdot e^{-\lambda_\xi t}|^{2}\cdot\|e_\xi\|_{\mathcal H}^{2}}{|\lambda_\xi|^{2}\cdot|1-e^{-\lambda_\xi T}|^{2} }\\&
\leq C\sum_{\xi\in \mathcal I}{|(\mathcal L\varphi,e_\xi)_{\mathcal H}|+|(\mathcal L\psi,e_\xi)_{\mathcal H}|}<\infty, \,\,\,C=const>0.\end{align*}

From this and $\varphi, \psi\in \mathcal H^{1}$, we obtain
\begin{align*}\|u(x,t)\|_{C^1([0,T],\mathcal H)}&=\max_{t\in [0,T]}\|u(x,t)\|_{\mathcal H}+\max_{t\in [0,T]}\|u_t(x,t)\|_{\mathcal H}<\infty.\end{align*}

Similarly for $f(x)$, we obtain the estimate
$$
\|f\|_{\mathcal H}^{2}\leq C\|\mathcal{L}\varphi\|_{\mathcal H}^{2}+C\sum_{\xi\in \mathcal I}{|(\mathcal L\varphi,e_\xi)_{\mathcal H}|^{2}+|(\mathcal L\psi,e_\xi)_{\mathcal H}|^{2}}<\infty.
$$
Existence of the solution of Problem \ref{PR-01} is proved.

\subsection{Proof of the uniqueness result}
Suppose that there are two solutions $\{u_1(x,t),f_1(x)\}$ and $\{u_2(x,t),f_2(x)\}$ of Problem \ref{PR-01}.
Denote $$u(x,t)=u_1(x,t)-u_2(x,t)$$ and
$$f(x)=f_1(x)-f_2(x).$$ Then the functions $u(x,t)$ and $f(x)$ satisfy Equation \eqref{1} and homogeneous conditions \eqref{2} and \eqref{3}.
We also have
 \begin{equation}\label{6} u_\xi(t)=(u(x,t),e_\xi(x))_{\mathcal H},\,\,\, \xi\in\mathcal I,\end{equation} and
\begin{equation}\label{7} f_\xi=(f(x),e_\xi(x))_{\mathcal H},\,\,\, \xi\in\mathcal I.\end{equation}
Applying the operator $\frac{d}{dt}$ to \eqref{6}, we have $$\frac{du_\xi(t)}{dt}=(u_t(x,t),e_\xi(x))_{\mathcal H}=(\mathcal{L}u(x,t)+f(x),e_\xi(x))_{\mathcal H}=\lambda_\xi u_\xi(t)+f_\xi.$$

Thus we get \begin{equation}\frac{du_\xi(t)}{dt}=\lambda_\xi u_\xi(t)+f_\xi,\end{equation} an ordinary first-order differential equation. The general solution of this equation is:

$$u_\xi(t)=A_\xi e^{\lambda_\xi t}-\frac{f_\xi}{\lambda_\xi},$$ where $A_\xi$ and $f_\xi$ are unknown constants. Using the homogeneous conditions \eqref{2} and \eqref{3} we obtain following conditions: $$u_\xi(0)=u_\xi(T)=0.$$ Using this we can find unknown constants $A_\xi$ and $f_\xi$. We first find  $A_\xi$: $$u_\xi(0)=A_\xi-\frac{f_\xi}{\lambda_\xi}=0,\,\,\, \textrm{then} \,\,\,A_\xi=\frac{f_\xi}{\lambda_\xi}.$$ Similarly, from
$$u_\xi(T)=A_\xi e^{\lambda_\xi T}-\frac{f_\xi}{\lambda_\xi}=0$$ we obtain $$\frac{f_\xi}{\lambda_\xi}(e^{\lambda_\xi T}-1)=0,$$ this implies $$f_\xi=0, A_\xi=0 \Longrightarrow u_\xi(t)\equiv0.$$
Further, by the completeness of the system $e_\xi$ in $\mathcal H$, we obtain $f(x)\equiv0, u(x,t)\equiv0.$
Uniqueness of the solution of the Problem \ref{PR-01} is proved.

\section{Inverse problem for the time-fractional diffusion equation}\label{InvFr}

The section deals with an inverse problem concerning the time-fractional diffusion equation.

\subsection{Preliminaries}\label{S-F}

Now, to formulate the problem, we need to define fractional differentiation operators.

\begin{defi}[\cite{Kilbas}] The left and right Riemann--Liouville
fractional integrals $I_{a+} ^\alpha$ and $I_{b-} ^\alpha$ of order $\alpha\in\mathbb R$ ($\alpha>0$) are given by
$$
I_{a+} ^\alpha  \left[ f \right]\left( t \right) = {\rm{
}}\frac{1}{{\Gamma \left( \alpha \right)}}\int\limits_a^t {\left(
{t - s} \right)^{\alpha  - 1} f\left( s \right)} ds, \,\,\, t\in(a,b],
$$
and
$$ I_{b-}
^\alpha  \left[ f \right]\left( t \right) = {\rm{
}}\frac{1}{{\Gamma \left( \alpha \right)}}\int\limits_t^b {\left(
{s - t} \right)^{\alpha  - 1} f\left( s \right)} ds, \,\,\, t\in[a,b),
$$
respectively. Here $\Gamma$ denotes the Euler gamma function.
\end{defi}

\begin{defi}[\cite{Kilbas}] The left Riemann--Liouville
fractional derivative $D_{a+} ^\alpha$ of order $\alpha\in\mathbb R$ ($0<\alpha<1$) is defined by
$$
D_{a+} ^\alpha \left[ f \right]\left( t \right) = \frac{{d }}{{dt
}}I_{a+} ^{1 - \alpha } \left[ f \right]\left( t \right), \,\,\, \forall t\in(a, b].
$$
Similarly, the right Riemann--Liouville
fractional derivative $D_{b-} ^\alpha$ of order $\alpha\in\mathbb R$ ($0<\alpha<1$) is given by
$$
D_{b-}
^\alpha \left[ f \right]\left( t \right) = -
\frac{{d }}{{dt }}I_{b-} ^{1 - \alpha } \left[ f \right]\left( t
\right), \,\,\, \forall t\in[a, b).
$$
\end{defi}

\begin{defi}[\cite{Kilbas}] The left and right Caputo fractional derivatives of order $\alpha\in\mathbb R$ ($0<\alpha<1$) are defined by
$$
\mathcal{D}_{a+} ^\alpha  \left[ f \right]\left( t \right) = D_{a+}
^\alpha  \left[ f\left( t \right) - f\left( a \right) \right], \,\,\, t\in(a, b],
$$
and
$$
\mathcal{D}_{b-} ^\alpha
\left[ f \right]\left( t \right) = D_{b-} ^\alpha  \left[ f\left( t
\right) - f\left( b \right) \right], \,\,\, t\in[a, b),
$$
respectively.
\end{defi}
Now, we are in a way to state our problem.

\begin{prm}\label{pr1} Find the couple of functions $\left(u\left( {x,t}
\right), f\left( x \right)\right)$ satisfying the equation
\begin{equation}\label{2.1}
\mathcal{D}_{0+,t}^{\alpha} u \left( {x,t} \right) + \mathcal{L} u(x, t) = f\left( x \right),
\end{equation} in the domain $(x,t)\in\Omega_t=\{x\in Q\subset \mathbb{R}^d,\,d\geq 1,\,\, t\in[0,T]\},$
under the conditions \begin{equation}\label{2.2}u\left(
{x,0} \right) = \varphi \left( x \right),\  x \in \bar Q,\end{equation} \begin{equation}\label{2.2*}u\left( {x,T} \right) =
\psi \left( x \right),\ x \in \bar Q,
\end{equation} where $\varphi \left( x \right)$ and $\psi
\left( x \right)$ are sufficiently smooth functions, $\mathcal L$ is a linear self-adjoint operator with a discrete spectrum.\end{prm}

If $\alpha =1,$ then equation \eqref{2.1} coincides with the classical heat equation. The heat equation also describes the diffusion process. So, the equation of the form \eqref{2.1} with fractional derivatives with respect to the time variable is called the sub-diffusion equation \cite{Uchaikin}. This equation describes the slow diffusion. When $\alpha =\frac{1}{2}$ the equation was interpreted by Nigmatullin \cite{Nigmatullin} within a percolation (pectinate) model. The solution (in an unbounded domain in the space variable) was investigated by Mainardi \cite{Mainardi} and others by means of integral transformations.

\begin{defi}[\cite{Carvalho}] Let $X$ be a Banach space. We say that $u\in C^\alpha([0,T],X)$ if $u\in C([0,T],X)$ and $\mathcal{D}^\alpha_tu\in C([0,T],X).$

Observe that, viewed as a subspace of $C([0,T], X),$ the space $C^\alpha([0, T], X)$ is a Banach space.
\end{defi}
A generalised solution of Problem \ref{pr1} is the pair of functions $(u(x,t),f(x)),$ where $u=u(x,t)\in C^\alpha([0,T],\mathcal H)\cap C([0,T],\mathcal H^{1})$ and $f=f(x)\in \mathcal H$.

For the considered Problem \ref{pr1}, the following theorem holds true.

\begin{thm}
\label{Th-01-F}
Let $\varphi,\psi\in \mathcal H^{1}.$ Then the generalised solution of the Problem \ref{pr1}, exists, is unique, and can be written in the form
$$
u(x,t)=\varphi(x)+\sum_{\xi\in \mathcal I}\frac{[(\varphi,e_\xi)_{\mathcal H}-(\psi,e_\xi)_{\mathcal H}](E_{\alpha,1}\left({-\lambda_\xi t^\alpha}\right)-1)e_\xi(x)}{(1-E_{\alpha,1}\left({-\lambda_\xi t^\alpha}\right)) },
$$

$$
f(x)=\mathcal L\varphi(x)-\sum_{\xi\in \mathcal I}\frac{[(\mathcal L\varphi,e_\xi)_{\mathcal H}-(\mathcal L\psi,e_\xi)_{\mathcal H}]e_\xi(x)}{(1-E_{\alpha,1}\left({-\lambda_\xi t^\alpha}\right))},
$$ where $E_{\alpha,1}(z)$ is the Mittag-Leffler function: $$E_{\alpha,1}(z)=\sum\limits_{m=0}^\infty\frac{z^m}{\Gamma(\alpha m+1)}.$$\end{thm} For different properties of the Mittag-Leffler function $E_{\alpha,1}(z)$ see e.g. \cite{Kilbas}.

\subsection{Proof of Theorem \ref{Th-01-F}}
We give the full proof for Problem \ref{PR-01}.
\subsubsection{Existence of solution} We want to find a generalised solution by Fourier method, we have the eigenvalues $\lambda_\xi$ and  eigenfunctions system $e_\xi(x)$  of operator $\mathcal L$ on the space $\mathcal H$. Eigenfunctions system $e_\xi(x)$ is an orthonormal basis in $\mathcal H$, so that the functions $u(x,t)$ and $f(x)$ can be expanded as follows:
\begin{equation}\label{6.1}u\left( {x,t} \right) = \sum\limits_{\xi\in \mathcal{I}}  {u_{\xi} \left( t \right)e_\xi(x)},\end{equation} \begin{equation}\label{6.2}f\left( x \right) =\sum\limits_{\xi\in \mathcal{I}} {f_{\xi}e_\xi(x)},\end{equation} where $u_{\xi} \left( t \right), f_{\xi}$ are unknown. Substituting (\ref{6.1}) and (\ref{6.2}) into (\ref{2.1}), we obtain the following equations for
the unknown functions $u_{\xi} \left( t \right)$ and the constants $f_{\xi}:$ $$\mathcal{D}^\alpha_{0+} u_{\xi} \left( t \right) + \lambda_{\xi} u_{\xi}\left( t \right) = f_{\xi},\,\,\, \xi\in\mathcal{I}.$$ By solving
this equation (see \cite{Kilbas}), we obtain  $$u_{\xi} \left( t \right) = \frac{{f_{\xi}}}{{\lambda_{\xi}}} + C_{\xi} E_{\alpha,1}\left({ - \lambda_{\xi} t^\alpha}\right),$$ where the constants $C_{\xi}$ and $f_{\xi}$ are unknown and $E_{\alpha,1}(z)$ is the Mittag-Leffler function \cite{Kilbas}: $$E_{\alpha,1}(z)=\sum\limits_{m=0}^\infty\frac{z^m}{\Gamma(\alpha m+1)}.$$ To find these constants, we use conditions (\ref{2.2}). Let $\varphi_\xi, \psi _\xi$ be the coefficients of the expansions of $\varphi \left( x \right)$ and $\psi \left( x \right)$
$$
\varphi_\xi  =  \left({\varphi, e_\xi}\right)_\mathcal{H},\,\,\psi_{\xi}  =  \left({\psi,e_\xi}\right)_\mathcal{H},
$$
where $(\cdot,\cdot)_\mathcal{H}$ is the inner product of the Hilbert space $\mathcal{H}.$

Then for $C_{\xi}$ we have
$$u_{\xi} \left( 0 \right) = \frac{{f_{\xi}}}{\lambda_{\xi}} + C_{\xi}  = \varphi_{\xi},
$$
$$u_\xi \left( T \right) =\frac{{f_{\xi}}}{\lambda_{\xi}} + C_{\xi} E_{\alpha,1}\left({ - \lambda_{\xi} T^\alpha}\right)  = \psi _\xi,
$$
$$\varphi _{\xi}  - C_{\xi}  + C_{\xi} E_{\alpha,1}\left({ - \lambda_{\xi} T^\alpha} \right) = \psi _{\xi}.
$$
Thus, we get
$$C_{\xi}  = \frac{{\varphi _{\xi}  - \psi _{\xi} }}{{1 - E_{\alpha,1}\left({ - \lambda_{\xi} T^\alpha}\right)}}.
$$
The unknowns $f_{\xi}$ can be represented as
$$f_\xi  = \lambda_{\xi} (\varphi_{\xi}  - C_{\xi}).
$$

Substituting $u_{\xi}\left( t \right),$ $f_{\xi}$ into (\ref{6.1}) and (\ref{6.2}), we find
$$u\left( {x,t} \right)=\varphi \left( x \right) + \sum\limits_{\xi\in \mathcal{I}}  {C_{\xi}\left( {E_{\alpha,1}\left({ - \lambda_{\xi} t^\alpha}\right)  - 1} \right)e_{\xi}}(x).$$

Using the self-adjointness property of the operator $\mathcal{L},$ we have $$(\mathcal{L}\varphi,e_\xi)_{\mathcal H}=(\varphi,\mathcal{L}e_\xi)_{\mathcal H}$$ we know that $\mathcal{L}e_\xi=\lambda_\xi e_\xi$ using this we obtain
$$(\varphi,e_\xi)_{\mathcal H}=\frac{(\mathcal{L}\varphi,e_\xi)_{\mathcal H}}{\lambda_\xi},$$ and for $\psi(x)$ we can write it in a similar way. Substituting these equalities into formula of $C_\xi$ we can get that $$C_\xi=\frac{(\mathcal L\varphi,e_\xi)_{\mathcal H}-(\mathcal L\psi,e_\xi)_{\mathcal H}}{\lambda_\xi(1-E_{\alpha,1}{(-\lambda_\xi T^\alpha)})}.$$
Then $$u(x,t)=\varphi(x)+\sum_{\xi\in \mathcal I}\frac{[(\mathcal L\varphi,e_\xi)_{\mathcal H}-(\mathcal L\psi,e_\xi)_{\mathcal H}](E_{\alpha,1}{(-\lambda_\xi t^\alpha)}-1)e_\xi(x)}{\lambda_\xi(1-E_{\alpha,1}{(-\lambda_\xi T^\alpha)}) }.$$
Similarly, $$f(x)=\mathcal L\varphi(x)-\sum_{\xi\in \mathcal I}\frac{[(\mathcal L\varphi,e_\xi)_{\mathcal H}-(\mathcal L\psi,e_\xi)_{\mathcal H}]e_\xi(x)}{(1-E_{\alpha,1}{(-\lambda_\xi T^\alpha)})}.$$

The following Mittag-Leffler function's estimate is known \cite{Simon}:
$$\frac{1}{1+\Gamma(1-\alpha)z}\leq E_{\alpha,1}(-z)\leq \frac{1}{1+\Gamma(1-\alpha)^{-1}z},\,\,z>0.$$
From this inequality it follows that
$$0< E_{\alpha,1}(-z)< 1,\,\,z>0.$$
Now, for the convergence of the series, using $\|e_\xi\|_{\mathcal H}=1$ and $\varphi, \psi\in \mathcal H^{1}$, we have the following estimate
\begin{align*}
\max_{t\in [0,T]}\|u(x,t)\|_{\mathcal H}^{2}\leq C\|\varphi(x)\|_{\mathcal H}^{2}+C\sum_{\xi\in \mathcal I}|(\varphi,e_\xi)_{\mathcal H}|^{2}+|(\psi,e_\xi)_{\mathcal H}|^{2}<\infty.
\end{align*} We also have
\begin{align*}
\max_{t\in [0,T]}\|\mathcal{D}^\alpha_{0+,t}u(x,t)\|_{\mathcal H}^{2}&\leq\max_{t\in [0,T]}\sum_{\xi\in \mathcal I}\frac{|(\mathcal L\varphi,e_\xi)_{\mathcal H}-(\mathcal L\psi,e_\xi)_{\mathcal H}|^{2}\cdot|-\lambda_\xi\cdot E_{\alpha,1}{(-\lambda_\xi t^\alpha)}|^{2}\cdot\|e_\xi\|_{\mathcal H}^{2}}{\lambda_\xi^{2}\cdot|1-E_{\alpha,1}{(-\lambda_\xi T^\alpha)}|^{2} }\\&
\leq C\sum_{\xi\in \mathcal I}{|(\mathcal L\varphi,e_\xi)_{\mathcal H}-(\mathcal L\psi,e_\xi)_{\mathcal H}|^{2}}\\&\leq C\sum_{\xi\in \mathcal I}{|(\mathcal L\varphi,e_\xi)_{\mathcal H}|^{2}+|(\mathcal L\psi,e_\xi)_{\mathcal H}|^{2}}<\infty.
\end{align*}

From this and $\varphi, \psi\in \mathcal H^{1}$, we obtain
\begin{align*}\|u(x,t)\|_{C^\alpha([0,T],\mathcal H)}^{2}&=\max_{t\in [0,T]}\|u(x,t)\|_{\mathcal H}^{2}+\max_{t\in [0,T]}\|\mathcal{D}^\alpha_{0+,t}u(x,t)\|_{\mathcal H}^{2}\\&\leq C\|\varphi(x)\|_{\mathcal H}^{2}+C\sum_{\xi\in \mathcal I}|(\varphi,e_\xi)_{\mathcal H}|^{2}+|(\psi,e_\xi)_{\mathcal H}|^{2}
\\&+C\sum_{\xi\in \mathcal I}{|(\mathcal L\varphi,e_\xi)_{\mathcal H}|^{2}+|(\mathcal L\psi,e_\xi)_{\mathcal H}}|^{2}<\infty.
\end{align*}

Similarly for $f(x)$, we obtain the estimate
$$
\|f\|_{\mathcal H}^{2}\leq C\|\mathcal{L}\varphi\|_{\mathcal H}^{2}+C\sum_{\xi\in \mathcal I}{|(\mathcal L\varphi,e_\xi)_{\mathcal H}|^{2}+|(\mathcal L\psi,e_\xi)_{\mathcal H}}|^{2}<\infty.
$$
Existence of the solution of the Problem \ref{pr1} is proved.

\subsubsection{Uniqueness of solution}
Hence the obtained solution satisfies the equation (\ref{2.1}) point-wise; by construction, it satisfies the conditions (\ref{2.2})-(\ref{2.2*}).

Suppose that there
are two solutions $\left\{ {u_1 \left( {x,t} \right),f_1 \left( x
\right)} \right\}$ and $\left\{ {u_2 \left( {x,t} \right), f_2
\left( x \right)} \right\}$ of Problem \ref{pr1}. Denote $$u\left( {x,t}
\right) = u_1 \left( {x,t} \right) - u_2 \left( {x,t} \right)$$
and $$f\left( x \right) = f_1 \left( x \right) - f_2 \left( x
\right).$$ Then the functions $u\left( {x,t} \right)$ and $f\left(
x \right)$ satisfy (\ref{2.1}) and homogeneous conditions
(\ref{2.2}) and (\ref{2.2*}).

Let \begin{equation}\label{5.2}u_{\xi} \left(t \right) = \left( u(x,t),e_\xi(x)\right)_\mathcal{H}, \;\xi \in \mathcal{H},\end{equation} and
\begin{equation}\label{5.5} f_{\xi} = \left(f(x),e_\xi(x)\right)_\mathcal{H}, \;\xi \in \mathcal{I}.\end{equation}
Applying the operator $\mathcal{D}^\alpha_{0+}$ to the equation (\ref{5.2}), we have \begin{align*}\mathcal{D}^\alpha_{0+} u_{\xi} \left( t\right) = \left(\mathcal{D}^\alpha_{0+,t}u(x,t),e_\xi(x)\right)_\mathcal{H} =- \left\langle \mathcal{L}u(x,t),e_\xi(x)\right\rangle_\mathcal{H} + f_\xi.\end{align*}

By self-adjointness and taking into account the homogeneous
conditions (\ref{2.2}) and (\ref{2.2*}), we obtain $$\mathcal{D}^\alpha_{0+} u_{\xi} \left( t\right)+\lambda_\xi u_\xi(t) = f_{\xi}\ ,\ u_{\xi}\left( 0 \right) = 0\ ,\ u_{\xi}\left( T \right) =0.$$

Consequently, $f_{\xi} \equiv 0, u_{\xi} \left( t \right) \equiv 0.$

Further, by the completeness of the system $\left\{e_\xi(x)\right\}_{\xi\in\mathcal{I}}$ in $\mathcal{H}$ we obtain $$f\left( t \right) \equiv 0,u\left(
{x,t} \right) \equiv 0,0 \le t \le T, x\in \bar Q.$$

Hence, uniqueness of the solution of Problem \ref{pr1} is proved.

\section{Inverse time-fractional diffusion problem for the hypoelliptic operators}

In this section we show how the obtained results can be applied for the analysis of the corresponding inverse problems for a large class of hypoelliptic diffusions and subdiffusions.

\subsection{Graded Lie groups}
We start by recalling following Folland and Stein \cite{FS-book} or \cite[Section 3.1]{FR16} and give some definitions and notations.
A Lie algebra $\mathfrak g$ is graded if it is endowed with a vector space decomposition
$$
\mathfrak g=\bigoplus_{j=1}^\infty V_j\;\textrm{ such that } [V_i,V_j]\subset V_{i+j},
$$
where all but finitely many of $V_j$'s are zero.
Consequently, we call that a connected simply connected Lie group $\G$ is graded if its Lie algebra $\mathfrak g$ is graded.
When the first stratum $V_1$ generates $\mathfrak g$ as an algebra, we get a stratified case of $\G$.

Graded Lie groups are homogeneous Lie groups with dilations. Define the operator $A$ by setting
$AX=\nu_j X$ for $X\in V_j$. Then the dilations on $\mathfrak g$ are defined by
$$
D_r:={\rm Exp}(A\ln r), \; r>0.
$$
The homogeneous dimension $Q$ of the graded Lie group $\G$ is defined by
$$
Q:=\nu_1+\ldots+\nu_n={\rm Tr}\, A.
$$

In what follows, we assume that $\G$ is a graded Lie group. Rockland operators are firstly defined in \cite{Rockland} through the representations. By following \cite[Definition 4.1.1]{FR16},
we call that $\R$ is a Rockland operator on the graded Lie group $\G$ if $\R$ is a left-invariant differential operator which is homogeneous of a positive order $\nu\in\mathbb N$ and satisfies the Rockland condition: for all representations $\pi\in \widehat{\G}$, excluding the trivial one, $\pi(\R)$ is injective on $\mathcal H^{\infty}_{\pi}$, i.e., from
$$
\pi(\R)v=0
$$
it follows that $v=0$ for arbitrary $v\in \mathcal H^{\infty}_{\pi}$.

We denote by $\widehat{\G}$ the unitary dual of $\G$, $\mathcal H^{\infty}_{\pi}$ is the space of smooth vectors of the representation $\pi\in\widehat{\G}$, and $\pi(\R)$ is the infinitesimal representation of $\R$,
see \cite[Definition 1.7.2]{FR16}.
For more information on graded Lie groups and Rockland operators the readers are referred to \cite[Chapter 4]{FR16}.

Let $\pi$ be a representation of the graded Lie group $\G$ on the separable Hilbert space $\mathcal{H}_{\pi}$. We say that $v\in \mathcal{H}_{\pi}$ is smooth if the function
$$
\G\ni x \mapsto \pi(x)v\in \mathcal{H}_{\pi}
$$
is of class $C^{\infty}$. We denote by $\mathcal{H}_{\pi}^{\infty}$ the space of all smooth vectors of a representation $\pi$.

In the paper \cite{HJL85} Hulanicki, Jenkins and Ludwig showed that the spectrum of $\pi(\R)$ is purely discrete and positive. Hence we can choose an orthonormal basis for $\mathcal{H}_{\pi}$ such that the infinite matrix related to the (self-adjoint) operator $\pi(\R)$ has the following form
\begin{equation}\label{pi_R_matrix}
\pi(\R)=\begin{pmatrix} \pi_{1}^{2} & 0 & \ldots &\ldots \\
0 & \pi_{2}^{2} & 0 & \ldots\\
\vdots & 0 & \ddots &\\
\vdots & \vdots &  & \ddots \end{pmatrix},
\end{equation}
where $\pi\in \widehat{\G}\backslash \{1\}$ and $\pi_{j}\in \mathbb{R}_{>0}$.

For $f\in L^{1}(\G)$ and $\pi \in \widehat{G}$, let us define the group Fourier transform of $f$ at $\pi$ by
$$
\mathcal{F}_{\G}f(\pi)\equiv \widehat{f}(\pi)\equiv \pi(f):=\int_{\G}f(x)\pi(x)^{*}dx,
$$
with the integration against the biinvariant Haar measure on the graded Lie group $\G$, which implies that a linear mapping $\widehat{f}(\pi)$ from $\mathcal{H}_{\pi}$ to itself can be represented by an infinite matrix once we choose a basis for $\mathcal{H}_{\pi}$. Consequently, we obtain
$$
\mathcal{F}_{\G}(\R f)(\pi)=\pi(\R)\widehat{f}(\pi).
$$
From now on, when we write $\widehat{f}(\pi)_{m,k}$, we will be using the same basis in $\mathcal{H}_{\pi}$ as the one giving \eqref{pi_R_matrix}.

In \cite{CG90} by Kirillov's orbit method, the authors showed that the Plancherel measure $\mu$ on $\widehat{\G}$ can be constructed explicitly. This means that we can have the Fourier inversion formula. Furthermore, $\pi(f)=\widehat{f}(\pi)$ is the Hilbert-Schmidt operator, i.e.
$$
\|\pi(f)\|^{2}_{{\rm HS}}={\rm Tr}(\pi(f)\pi(f)^{*})<\infty,
$$
and $\widehat{\G}\ni \pi \mapsto \|\pi(f)\|^{2}_{{\rm HS}}$ is an integrable function with respect to $\mu$. Moreover, the Plancherel formula holds (see e.g. \cite{CG90} or \cite{FR16}):
\begin{equation}\label{planch_for}
\int_{\G}|f(x)|^{2}dx=\int_{\widehat{\G}}\|\pi(f)\|^{2}_{{\rm HS}}d\mu(\pi).
\end{equation}

In \cite{HN-79} Helffer and Nourrigat showed that a left-invariant differential operator $\R$ of homogeneous positive degree $\nu\in\mathbb N$ satisfies the Rockland condition if and only if it is hypoelliptic. Such operators we call Rockland operators.

The  Sobolev spaces $H^s_\R(\G)$, $s\in\mathbb R$, associated to
positive Rockland operators $\R$ have been analysed in \cite{FR:Sobolev} (see also \cite{FR16}). One of the definitions of Sobolev spaces is
\begin{equation*}\label{EQ:HsL-00-g}
H^s(\G):=H^s_\R(\G):=\left\{ f\in\Dcal'(\G): (I+\R)^{s/\nu}f\in
L^2(\G)\right\},
\end{equation*}
with the norm $\|f\|_{H^s_\R(\G)}:=\|(I+\R)^{s/\nu}f\|_{L^2(\G)},$ for a positive Rockland operator of homogeneous degree $\nu$. It is known that these spaces do not depend on a particular choice of a Rockland operator used in the above definition.
We refer to \cite{FR16} for details of the Fourier analysis on graded Lie groups.

Now we state the main problem of this section.

\begin{prm}\label{R: Pr-1} Let $\G$ be a graded Lie group of homogeneous dimension $Q\geq 3$ and let $\R$ be a positive self--adjoint Rockland operator acting on $L^{2}(\G)$. In a domain $\Omega=\{(t,x):\,(0,T)\times \G\}$ we seek a couple of functions $\left( u(t), f \right)$ satisfying the equation
\begin{equation}\label{R: 2.1}
\mathcal{D}_{0+,t}^{\alpha} u \left( {x,t} \right) + \R u(x, t) = f\left( x \right),
\end{equation}
under the conditions
\begin{equation}\label{R: 2.2}
u\left({x,0} \right) = \varphi \left( x \right),\  x \in \G,
\end{equation}
\begin{equation}\label{R: 2.2*}
u\left( {x,T} \right) =
\psi \left( x \right),\ x \in \G,
\end{equation}
where $\mathcal{D}_{0+,t}^{\alpha}$ is a Caputo fractional derivative, $0 < \alpha \leq 1$.

We seek a solution $u\in C([0,T]; L^{2}(\G))$ of the problem \eqref{R: 2.1}--\eqref{R: 2.2*} such that $\mathcal{D}_{0+,t}^{\alpha} u\in C([0,T]; L^{2}(\G))$, $\R u\in C([0,T]; L^{2}(\G))$, and $f\in C([0,T]; L^{2}(\G))$.
\end{prm}

\begin{thm}
\label{R: Th-01-F}
Assume that $\R$ is a positive Rockland operator of homogeneous order $\nu$. Let $\varphi,\psi\in H^{\nu}(\G).$ Then there exists a unique solution $(u, f)$ such that $f\in L^2(\G)$, $u\in C([0,T]; H^{\nu}(\G))$, and $\mathcal{D}_{0+,t}^{\alpha} u\in C([0,T]; L^2(\G))$ of Problem \ref{R: Pr-1}, and can be written in the form
\begin{equation}
\label{R: E2}
u(t,x)=\int_{\widehat{\G}}\mathrm{Tr}[\widehat{K}(t, \pi)\pi(x)]d\mu(\pi),
\end{equation}
and
$$
f(x)=\int_{\widehat{\G}}\mathrm{Tr}[\widehat{F}(\pi)\pi(x)]d\mu(\pi),
$$
where
$$
\widehat{K}(t, \pi)_{l, k}=\widehat {\varphi}(\pi)_{l,k}+\frac{[\widehat {\varphi}(\pi)_{l,k}-\widehat {\psi}(\pi)_{l,k}](E_{\alpha,1}\left({-\pi_{l}^{2} t^\alpha}\right)-1)}{(1-E_{\alpha,1}\left({-\pi_{l}^{2} T^\alpha}\right))},
$$
and
$$
\widehat{F}(\pi)_{l, k}=\pi_{l}^{2}\widehat {\varphi}(\pi)_{l,k} - \frac{\pi_{l}^{2}[\widehat {\varphi}(\pi)_{l,k}-\widehat {\psi}(\pi)_{l,k}]}{1-E_{\alpha,1}\left({-\pi_{l}^{2} T^\alpha}\right)},
$$
for all $\pi \in \widehat{\G}$ and $l,k\in \mathbb{N}$, where $E_{\alpha,1}(z)$ is the Mittag-Leffler function.
\end{thm}

\subsection{Proof of Theorem \ref{R: Th-01-F}}
\subsubsection{Proof of the existence result.}
We give a full proof of Problem \ref{R: Pr-1}.
Let us take the group Fourier
transform of \eqref{R: 2.1} with respect to $x\in\G$ for all $\pi\in\widehat{G}$, that is,
\begin{equation}\label{fou_cauchy1}
\mathcal{D}_{0+,t}^{\alpha} \widehat {u}(t,\pi)+ \pi(\R)\widehat{u}(t,\pi)=\widehat{f}(\pi).
\end{equation}
Taking into account \eqref{pi_R_matrix}, we rewrite the matrix equation \eqref{fou_cauchy1} componentwise as an infinite system of equations of the form
\begin{equation}\label{R: 6}
\mathcal{D}_{0+,t}^{\alpha}\widehat {u}(t,\pi)_{l,k}+\pi_{l}^{2}\widehat{u}(t,\pi)_{l,k}=\widehat{f}(\pi)_{l,k},
\end{equation}
for all $\pi \in \widehat{\G}$, and any $l, k\in \mathbb{N}$. Now let us decouple the system given by the matrix equation \eqref{fou_cauchy1}. For this, we fix an arbitrary representation $\pi$, and a general entry $(l, k)$ and we treat each equation given by \eqref{R: 6} individually.

According to \cite{LG99}, the solutions of the equations \eqref{R: 6} satisfying initial conditions
\begin{equation}
\label{R: 6*}
\widehat {u}(0,\pi)_{l,k}=\widehat{\varphi}(\pi)_{l,k}, \, \widehat {u}(T,\pi)_{l,k}=\widehat{\psi}(\pi)_{l,k},
\end{equation}
can be represented in the form
\begin{equation}
\label{R: 7}
\widehat {u}(t,\pi)_{l,k} = \widehat {\varphi}(\pi)_{l,k}+\frac{[\widehat {\varphi}(\pi)_{l,k}-\widehat {\psi}(\pi)_{l,k}](E_{\alpha,1}\left({-\pi_{l}^{2} t^\alpha}\right)-1)}{(1-E_{\alpha,1}\left({-\pi_{l}^{2} T^\alpha}\right))},
\end{equation}
$$
\widehat{f}(\pi)_{l,k}=\pi_{l}^{2}\widehat {\varphi}(\pi)_{l,k} - \frac{\pi_{l}^{2}[\widehat {\varphi}(\pi)_{l,k}-\widehat {\psi}(\pi)_{l,k}]}{1-E_{\alpha,1}\left({-\pi_{l}^{2} T^\alpha}\right)},
$$
for all $\pi \in \widehat{\G}$ and any $l,k\in \mathbb{N}$, where $E_{\alpha,1}(z)$ is the Mittag-Leffler function \cite{Kilbas}:
$$
E_{\alpha,1}(z)=\sum\limits_{m=0}^\infty\frac{z^m}{\Gamma(\alpha m+1)}.
$$

Then there exists a solution of Problem \ref{R: Pr-1}, and it can be written as
\begin{equation}
\label{R: E2}
u(t,x)=\int_{\widehat{\G}}\mathrm{Tr}[\widehat{K}(t, \pi)\pi(x)]d\mu(\pi),
\end{equation}
$$
f(x)=\int_{\widehat{\G}}\mathrm{Tr}[\widehat{F}(\pi)\pi(x)]d\mu(\pi),
$$
where
$$
\widehat{K}(t, \pi)_{l, k}=\widehat {\varphi}(\pi)_{l,k}+\frac{[\widehat {\varphi}(\pi)_{l,k}-\widehat {\psi}(\pi)_{l,k}](E_{\alpha,1}\left({-\pi_{l}^{2} t^\alpha}\right)-1)}{(1-E_{\alpha,1}\left({-\pi_{l}^{2} T^\alpha}\right))},
$$
$$
\widehat{F}(\pi)_{l, k}=\pi_{l}^{2}\widehat {\varphi}(\pi)_{l,k} - \frac{\pi_{l}^{2}[\widehat {\varphi}(\pi)_{l,k}-\widehat {\psi}(\pi)_{l,k}]}{1-E_{\alpha,1}\left({-\pi_{l}^{2} T^\alpha}\right)},
$$
for all $\pi \in \widehat{\G}$ and $l,k\in \mathbb{N}$.

We note that the above expressions ares well-defined in view of $\varphi,\psi\in \H^{\nu}(\G)$. Finally, based on \eqref{R: 7}, we rewrite our formal solution as \eqref{R: E2}.

\subsubsection{Convergence of the formal solution.}
Here, we prove convergence of the obtained infinite series corresponding to functions $u(x,t)$, $\mathcal{D}_{0+,t}^{\alpha} u(x,t)$, and $\R u(x,t)$.

Thus, since for any Hilbert-Schmidt operator $A$ one has
$$\|A\|^{2}_{{\rm HS}}=\sum_{l,k}|(A\phi_{l},\phi_{k})|^{2}$$
for any orthonormal basis $\{\phi_{1},\phi_{2},\ldots\}$, then we can consider the infinite sum over $l,k$ of the inequalities provided by \eqref{R: 7}, we have
\begin{equation}
\label{eq4}
\|\widehat {u}(t,\pi)\|^{2}_{{\rm HS}} \leq C (\|\widehat{\varphi}(\pi)\|^{2}_{{\rm HS}}+\|\widehat{\psi}(\pi)\|^{2}_{{\rm HS}}),
\end{equation}
\begin{equation}
\label{eq4}
\|\mathcal{D}_{0+,t}^{\alpha}\widehat {u}(t,\pi)\|^{2}_{{\rm HS}} \leq C (\|\pi(\R)\widehat{\varphi}(\pi)\|^{2}_{{\rm HS}}+\|\pi(\R)\widehat{\psi}(\pi)\|^{2}_{{\rm HS}}),
\end{equation}
and
\begin{equation}
\label{eq4}
\|\widehat {\R u}(t,\pi)\|^{2}_{{\rm HS}} \leq C (\|\pi(\R)\widehat{\varphi}(\pi)\|^{2}_{{\rm HS}}+\|\pi(\R)\widehat{\psi}(\pi)\|^{2}_{{\rm HS}}).
\end{equation}
Thus, integrating both sides of \eqref{eq4} against the Plancherel measure $\mu$ on $\widehat{\G}$, then using the Plancherel identity \eqref{planch_for} we obtain
\begin{align*}
\|u\|_{C([0,T]; L^{2}(\G))} \leq C (\|\varphi\|_{L^{2}(\G)}+\|\psi\|_{L^{2}(\G)}),
\end{align*}
\begin{align*}
\|\mathcal{D}_{0+,t}^{\alpha} u\|_{C([0,T]; L^{2}(\G))} \leq C (\|\R\varphi\|_{L^{2}(\G)}+\|\R\psi\|_{L^{2}(\G)}),
\end{align*}
and
\begin{align*}
\|\R u\|_{C([0,T]; L^{2}(\G))} \leq C (\|\R\varphi\|_{L^{2}(\G)}+\|\R\psi\|_{L^{2}(\G)}).
\end{align*}

Similarly for $f$, we obtain the estimate
\begin{align*}
\|f\|_{L^{2}(\G)} \leq C (\|\varphi\|_{\H^{\nu}(\G)}+\|\R\psi\|_{\H^{\nu}(\G)}).
\end{align*}

The uniqueness result can be proved in analogy to the previous arguments.

\section{Appendix: Non-harmonic analysis of operators with discrete spectrum}
\label{App 01}

In this section we recall the necessary elements of the global Fourier analysis that has been developed in \cite{RT16}, and its applications to the spectral properties of operators in \cite{DRT17}. The space $$\mathcal{H}_\mathcal{L}^\infty:=\textrm{Dom}\left(\mathcal{L}^\infty\right)$$ is called the space of test functions for $\mathcal{L}.$ Here we define
$$\textrm{Dom}\left(\mathcal{L}^\infty\right):=\bigcap\limits_{k=1}^\infty\textrm{Dom}\left(\mathcal{L}^k\right),$$ where $\textrm{Dom}\left(\mathcal{L}^k\right)$ is the domain of the operator $\mathcal{L}^k,$ in turn defined as
$$\textrm{Dom}\left(\mathcal{L}^k\right):=\left\{f\in \mathcal{H}:\, \mathcal{L}^jf\in \textrm{Dom}\left(\mathcal{L}\right),\,j=0,1,2,...,k-1\right\}.$$
The Fr\'{e}chet topology of $\mathcal{H}^\infty_\mathcal{L}$ is given by the family of semi-norms
\begin{equation}\label{FA-1}\|\varphi\|_{\mathcal{H}^\infty_\mathcal{L}}:=\max_{j\geq k}\left\|\mathcal{L}^j\varphi\right\|_{\mathcal{H}},\, k\in \mathbb{N}_0,\, \varphi\in \mathcal{H}^\infty_\mathcal{L}.\end{equation}
Analogously to the operator $\mathcal{L}^*$ ($\mathcal{H}$-conjugate to $\mathcal{L}$), we introduce the space
$$\mathcal{H}_{\mathcal{L^*}}^\infty:=\textrm{Dom}\left(\left(\mathcal{L^*}\right)^\infty\right)$$
of test functions for $\mathcal{L^*},$ and we define
$$\textrm{Dom}\left(\left(\mathcal{L^*}\right)^\infty\right):=\bigcap\limits_{k=1}^\infty\textrm{Dom}\left(\left(\mathcal{L^*}\right)^k\right),$$
where $\textrm{Dom}\left(\left(\mathcal{L^*}\right)^k\right)$ is the domain of the operator $\left(\mathcal{L^*}\right)^k,$ in turn defined as
$$\textrm{Dom}\left(\left(\mathcal{L^*}\right)^k\right):=\left\{f\in \mathcal{H}:\, \left(\mathcal{L^*}\right)^jf\in \textrm{Dom}\left(\mathcal{L}^*\right),\,j=0,1,2,...,k-1\right\}.$$
The Fr\'{e}chet topology of $\mathcal{H}^\infty_{\mathcal{L}^*}$ is given by the family of semi-norms
\begin{equation}\label{FA-2}\|\varphi\|_{\mathcal{H}^\infty_{\mathcal{L}^*}}:=\max_{j\geq k}\left\|\left(\mathcal{L^*}\right)^j\varphi\right\|_{\mathcal{H}},\, k\in \mathbb{N}_0,\, \varphi\in \mathcal{C}^\infty_{\mathcal{L}^*}.\end{equation}
Now the space
$$\mathcal{H}^{-\infty}_{\mathcal{L}^*}:=\mathcal{L}\left(\mathcal{H}^{\infty}_{\mathcal{L}}, \mathbb{C}\right)$$
of linear continuous functionals on $\mathcal{H}^{\infty}_{\mathcal{L}}$ is called the space of $\mathcal{L}^*$-distributions. We can understand the continuity here in terms of the topology \eqref{FA-1}. For $w \in \mathcal{H}^{-\infty}_{\mathcal{L}^*}$ and $w \in \mathcal{H}^{\infty}_{\mathcal{L}},$ we shall also write
$$w(\varphi)=\langle w,\varphi\rangle.$$ For any $\psi\in \mathcal{H}^{\infty}_{\mathcal{L}^*},$ the functional $$\mathcal{H}^{\infty}_{\mathcal{L}}\ni \varphi \mapsto (\psi, \varphi)$$
is an ${\mathcal{L}^*}$-distribution, which gives an embedding $\psi\in \mathcal{H}^{\infty}_{\mathcal{L}^*}\hookrightarrow \mathcal{H}^{-\infty}_{\mathcal{L}^*}.$

Since the system of eigenfunctions $\{e_\xi: \xi \in \mathcal{I}\}$ of the operator $\mathcal{L}$ is a Riesz basis in $\mathcal{H}$ then its biorthogonal system $\{e^*_\xi: \xi\in \mathcal{I}\}$ is also a Riesz basis in $\mathcal H$ (see e.g. Bari \cite{B51}, as well as Gelfand \cite{G63}). Note that the function $e^*_\xi$ is an eigenfunction of the operator $\mathcal{L}^*$ corresponding to the eigenvalue $\bar{\lambda}_\xi$ for each $\xi\in \mathcal{I}.$ They satisfy the orthogonality relations
$$
(e_\xi, e^*_\eta) = \delta_{\xi, \eta},
$$
where $\delta_{\xi, \eta}$ is the Kronecker delta.

If $\mathcal{L}$ is self-adjoint, we clearly have $e^*_\xi=e_\xi.$

\section{Examples}
\label{SEC:examples}

Now as an illustration we give several examples of the settings where our inverse problems are applicable. Of course, there are many other examples, here we collect the ones for which different types of partial differential equations have particular importance.

\begin{itemize}
  \item {\bf Sturm-Liouville problem.}\\
  First, we describe the setting of the Sturm-Liouville operator. Let $l$ be an ordinary second order differential operator in $L^2(a,b)$ generated by the differential expression
\begin{equation}\label{SL} l(u)=-u''(x),\,\,a<x<b,\end{equation} and one of the boundary conditions \begin{equation}\label{SL_B} a_1u'(b)+b_1u(b)=0,\,a_2u'(a)+b_2u(a)=0,\end{equation} or \begin{equation}\label{SL_B1} u(a)=\pm u(b),\,u'(a)=\pm u'(b),\end{equation} where $a^2_1+a^2_2>0,\,b_1^2+b_2^2>0$ and $\alpha_j,\, \beta_j,\,j=1,2$ some real numbers.

It is known \cite{Naimark} that the Sturm-Liouville problem for equation \eqref{SL} with boundary conditions \eqref{SL_B} or with boundary conditions \eqref{SL_B1} is self-adjoint in $L^2(a,b).$ It is known that the self-adjoint problem has real eigenvalues and their eigenfunctions form a complete orthonormal basis in $L^2(a,b).$
\end{itemize}

\begin{itemize}
  \item {\bf Differential operator with involution.}\\
As a next example, we consider the differential operator with involution in $L^2(0,\pi)$ generated by the expression
\begin{equation}\label{DOI} l(u)=u''(x)-\varepsilon u''(\pi-x),\,\,0<x<\pi,\end{equation} and homogeneous Dirichlet conditions \begin{equation}\label{DOI_B} u(0)=0,\,u(\pi)=0,\end{equation} where $|\varepsilon|<1$ some real number.

The nonlocal functional-differential operator \eqref{DOI}-\eqref{DOI_B} is self-adjoint \cite{Kir1, torebek}. For $|\varepsilon|<1,$ the operator \eqref{DOI}-\eqref{DOI_B} has the following eigenvalues
\begin{align*}\lambda_{2k}=4(1+\varepsilon)k^2,\,k\in \mathbb{N}\,\,\, \textrm{and} \,\,\, \lambda_{2k+1}=(1-\varepsilon)(2k+1)^2,\,k\in \mathbb{N}\cup\{0\}\end{align*}
and corresponding eigenfunctions
\begin{align*}&u_{2k}(x)=\sqrt{\frac{2}{\pi}}\sin{2kx},\,k\in \mathbb{N},\\& u_{2k+1}(x)=\sqrt{\frac{2}{\pi}}\sin{(2k+1)x},\,k\in \mathbb{N}\cup \{0\}.\end{align*}
\end{itemize}

\begin{itemize}
  \item {\bf Fractional Sturm-Liouville operator.}\\
We consider the operator generated by the integro-differential expression
\begin{equation}\label{FSL}\ell (u)=\mathcal{D}_{a+}^\alpha D_{b-}^{\alpha}u,\,a<x<b,\end{equation}
and the conditions \begin{equation}\label{FSL_B}I_{b-}^{1-\alpha}u(a)=0,\,I_{b-}^{1-\alpha}u(b)=0,\end{equation} where
$\mathcal{D}_{a+}^\alpha$ is the left Caputo derivative of order $\alpha  \in
\left( {1/2,1} \right],$  $D_{b-}^\alpha$ is the right Riemann-Liouville derivative of order $\alpha \in
\left( {1/2,1} \right]$ and $I_{b-}^\alpha $ is the right Riemann-Liouville integral of order $\alpha \in
\left( {1/2,1} \right]$ (see. Subsection \ref{S-F} and \cite{Kilbas}).
The fractional Sturm-Liouville operator \eqref{FSL}-\eqref{FSL_B} is self-adjoint and positive in $L^2 (a, b)$ (see \cite{TT16}). The spectrum of fractional Sturm-Liouville operator \eqref{FSL}-\eqref{FSL_B} is
discrete, positive and real valued, and the system of eigenfunctions is a complete orthogonal basis in $L^2 (a, b).$ For more properties of the operator generated by the problem \eqref{FSL}-\eqref{FSL_B} we refer to \cite{TT18, TT18d, TT18a}.
\end{itemize}

\begin{itemize}
  \item {\bf Harmonic oscillator.}\\
For any dimension $d\geq1$, let us consider the harmonic oscillator,
$$
\L:=-\Delta+|x|^{2}, \,\,\, x\in\mathbb R^{d}.
$$
$\L$ is an essentially self-adjoint operator on $C_{0}^{\infty}(\mathbb R^{d})$. It has a discrete spectrum, consisting of the eigenvalues
$$
\lambda_{k}=\sum_{j=1}^{d}(2k_{j}+1), \,\,\, k=(k_{1}, \cdots, k_{d})\in\mathbb N^{d},
$$
and with the corresponding eigenfunctions
$$
\varphi_{k}(x)=\prod_{j=1}^{d}P_{k_{j}}(x_{j}){\rm e}^{-\frac{|x|^{2}}{2}},
$$
which are an orthogonal basis in $L^{2}(\mathbb R^{d})$. We denote by $P_{l}(\cdot)$ the $l$--th order Hermite polynomial, and
$$
P_{l}(\xi)=a_{l}{\rm e}^{\frac{|\xi|^{2}}{2}}\left(x-\frac{d}{d\xi}\right)^{l}{\rm e}^{-\frac{|\xi|^{2}}{2}},
$$
where $\xi\in\mathbb R$, and
$$
a_{l}=2^{-l/2}(l!)^{-1/2}\pi^{-1/4}.
$$
For more information, see for example \cite{NiRo:10}.
\end{itemize}

\begin{itemize}
  \item {\bf Anharmonic oscillator.}\\
Another class of examples -- anharmonic oscillators (see for instance
\cite{HR82}), operators on $L^2(\R)$ of the form
$$
\L:=-\frac{d^{2k}}{dx^{2k}} +x^{2l}+p(x), \,\,\, x\in\mathbb R,
$$
for integers $k,l\geq 1$ and with $p(x)$ being a polynomial of degree $\leq 2l-1$ with real coefficients. Other examples are a large class of harmonic and anharmonic oscillators in all dimension, see e.g. \cite{CDR18}.
\end{itemize}

\begin{itemize}
  \item {\bf Landau Hamiltonian in 2D.}\\
The next example is one of the simplest and most interesting models of the Quantum Mechanics, that is, the Landau Hamiltonian.

The Landau Hamiltonian in 2D is given by
\begin{equation} \label{eq:LandauHamiltonian}
\L:=\frac{1}{2} \p{\p{i\frac{\partial}{\partial x}-B
y}^{2}+\p{i\frac{\partial}{\partial y}+B x}^{2}},
\end{equation}
acting on the Hilbert space $L^{2}(\mathbb R^{2})$, where $B>0$ is some constant. The spectrum of
$\L$ consists of infinite number of eigenvalues (see \cite{F28, L30}) with
infinite multiplicity of the form
\begin{equation} \label{eq:HamiltonianEigenvalues}
\lambda_{n}=\p{2n+1}B, \,\,\, n=0, 1, 2, \dots \,,
\end{equation}
and the corresponding system of eigenfunctions (see \cite{ABGM15, HH13}) is
{\small
\begin{equation*}
\label{eq:HamiltonianBasis} \left\{
\begin{split}
e^{1}_{k, n}(x,y)&=\sqrt{\frac{n!}{(n-k)!}}B^{\frac{k+1}{2}}\exp\Big(-\frac{B(x^{2}+y^{2})}{2}\Big)(x+iy)^{k}L_{n}^{(k)}(B(x^{2}+y^{2})), \,\,\, 0\leq k, {}\\
e^{2}_{j, n}(x,y)&=\sqrt{\frac{j!}{(j+n)!}}B^{\frac{n-1}{2}}\exp\Big(-\frac{B(x^{2}+y^{2})}{2}\Big)(x-iy)^{n}L_{j}^{(n)}(B(x^{2}+y^{2})), \,\,\, 0\leq j,
\end{split}
\right.
\end{equation*}}
where $L_{n}^{(\alpha)}$ are the Laguerre polynomials given by
$$
L^{(\alpha)}_{n}(t)=\sum_{k=0}^{n}(-1)^{k}C_{n+\alpha}^{n-k}\frac{t^{k}}{k!}, \,\,\, \alpha>-1.
$$
Note that in \cite{RT17b, RT18, RT19a} the wave equation for the Landau Hamiltonian with a singular magnetic field was studied.

The reader is referred to \cite{MRT19a, MRT19b, RT17c, RT18c, RT19b} for more examples of operators $\L$ and its applications.
\end{itemize}

\begin{itemize}
  \item {\bf The restricted fractional Laplacian.}\\
  On the other hand, one can define a fractional Laplacian operator by using the integral representation in terms of hypersingular kernels already mentioned
  $$\left(-\Delta_{\mathbb{R}^n}\right)^sg(x)=C_{d,s}\, \textrm{P.V.} \int\limits_{\mathbb{R}^n}\frac{g(x)-g(\xi)}{|x-\xi|^{n+2s}}d\xi,$$ where $s\in (0,1).$

  In this case we materialize the zero Dirichlet condition by restricting the operator to act only on functions that are zero outside bounded domain $\Omega\subset\mathbb{R}^n.$ Caffarelli and Siro \cite{CS17} called the operator defined in such a way the restricted fractional Laplacian $\left(-\Delta_{\Omega}\right)^s.$ As defined, $\left(-\Delta_{\Omega}\right)^s$ is a self-adjoint operator on $L^2(\Omega),$ with a discrete spectrum $\lambda_{s,k}>0,\,\,k\in \mathbb{N}.$ The corresponding set of eigenfunctions $\{V_{s,k}(x)\}_{k\in \mathbb{N}},$ normalized in $L^2(\Omega).$
\end{itemize}

\begin{itemize}
  \item {\bf Hesienberg Harmonic and Anharmonic Oscillators.}\\
In \cite{RR18} the authors studied classes of operators yielding  harmonic and anharmonic oscillators on the Heisenberg group $\mathbf{H}_{n}$. These operators share one main feature: Every gradable Lie algebra equipped with a specific gradation admits a non-empty set of positive Rockland forms, consequently, one can associate to each positive Rockland form $P\in\mathfrak{u}(\mathfrak{h}_{n, 2})$ a family of anharmonic oscillators
$\{ \mathfrak{A}_{\mathbf{H}_{n}}^{k} \}_{k \in \mathbb{R} \setminus \{ 0 \}}$; the family of harmonic oscillators
$\{ \mathfrak{O}_{\mathbf{H}_{n}}^{k} \}_{k \in \mathbb{R} \setminus \{ 0 \}}$ was included as the special case in which the Dynin-Folland Lie algebra  $\mathfrak{h}_{n, 2}$ was equipped with the canonical homogeneous structure related to the natural stratification and $P=- (X_1^2 + \ldots X_{2n}^2 + Y_{2n+1}^2)$. Such operators have a discrete spectrum and also fall within the scope of this paper.
\end{itemize}

\section{Numerical illustrations}
\label{NumEx}
In this section we provide some numerical simulations. We deal with the case of operator \eqref{DOI}-\eqref{DOI_B}. Namely, in $L^{2}(0, \pi)$
consider the inverse source problem for the heat equation
\begin{equation}\label{CP-01}
u_t(x,t) - u''(x, t) + \varepsilon u''(\pi-x, t)=f(x), \,\,\, 0<x<\pi,
\end{equation}
with boundary conditions
\begin{equation}\label{BC-01}
u(0, t)=0, \, u(\pi, t)=0, \,\,\, t\in [0, T],
\end{equation}
and with the Cauchy problem
\begin{equation}\label{CP-02}
u(x,0)=\varphi(x), \,\,\, x\in [0, \pi],
\end{equation}
with some extra information
\begin{equation}\label{ED-01}
u(x,T)=\psi(x), \,\,\, x\in [0, \pi],
\end{equation}
where $-1<\varepsilon<1$ is some real number.

By Theorem \ref{Th-01}, from the Cauchy data $\varphi\in \mathcal H^2:=\dot{W}_{2}^{4}(0, \pi)$, that is, $\dot{W}_{2}^{4}(0, \pi):=\{\varphi, \varphi^{(4)}\in L^{2}(0, \pi) \,\,\, \hbox{with} \,\,\, \varphi(0)=\varphi(\pi)=\varphi''(0)=\varphi''(\pi)=0\},$ and from the observation function $\psi\in \dot{W}_{2}^{4}(0, \pi)$ we restore a unique solution $u\in C^2([0,T],L^{2}(0, \pi))\cap C([0,T], \dot{W}_{2}^{4}(0, \pi))$ and a source term $f\in L^{2}(0, \pi)$ by formulas
\begin{equation}\label{Sol-num}
u(x,t)= \varphi(x) + \sum_{k\in \mathbb N}\frac{\int_{0}^{\pi} l^{2}(\psi(s) - \varphi(s)) v_k(s) ds}{\lambda_k^2 (1-e^{-\lambda_k T}) } (1 -  e^{-\lambda_k t}) v_k(x),
\end{equation}
and
\begin{equation}\label{f-num}
f(x)= l(\varphi(x)) + \sum_{k\in \mathbb N}\frac{\int_{0}^{\pi} l^{2}(\psi(s) - \varphi(s)) v_k(s) ds  }{ \lambda_k (1-e^{-\lambda_k T})} v_k(x).
\end{equation}
Here, we denote
$$
l(\varphi(x)):= -\varphi''(x) + \varepsilon \varphi''(\pi-x).
$$

From Example \eqref{DOI}-\eqref{DOI_B}, we have a representation for the eigenvalues
\begin{align*}
\lambda_{k}=(1+(-1)^{k}\varepsilon)k^2, \, k\in \mathbb{N},
\end{align*}
and for the eigenfunctions
\begin{align*}
v_{k}(x)=\sqrt{\frac{2}{\pi}}\sin{kx}, \, k\in \mathbb{N}.
\end{align*}

\begin{figure}[ht!]
\begin{minipage}[h]{0.49\linewidth}
\center{\includegraphics[scale=0.4]{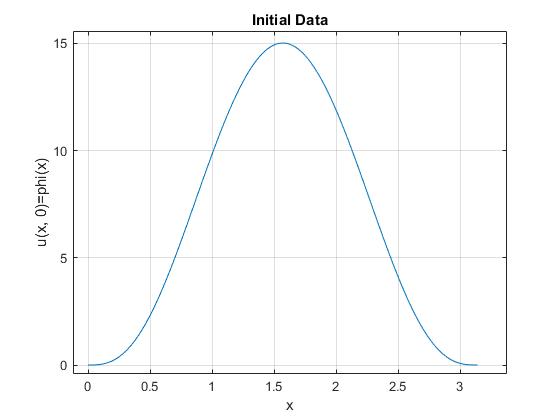}}
\end{minipage}
\caption{In this picture the temperature function at the initial point $0$ is drawn. The function is given by the formula \eqref{Test-IC}.} \label{fig1}
\end{figure}

In what follows, we consider and compare the cases $\varepsilon=0$ and $\varepsilon=0.9$. For numerical simulations, we choose
$$
\varphi(x)=x^{3}(\pi-x)^3, \,\,\, \psi=0,
$$
for all $x\in[0, \pi]$. Our problem is to cool (make $u(x, T)=0$ at some time $T$) a rod with the initial temperature
\begin{equation}
\label{Test-IC}
u(x, 0)=x^{3}(\pi-x)^3, \,\,\, x\in[0, \pi].
\end{equation}
Now, a question stands as follows: how should we choose a cooling function (a source term) $f$ to succeed the goal?

\begin{figure}[ht!]
\begin{minipage}[h]{0.49\linewidth}
\center{\includegraphics[scale=0.4]{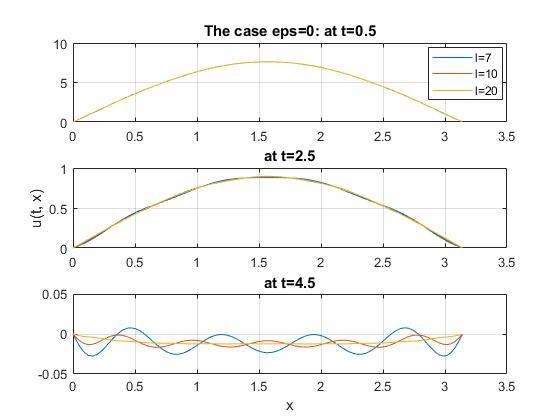}}
\end{minipage}
\hfill
\begin{minipage}[h]{0.49\linewidth}
\center{\includegraphics[scale=0.4]{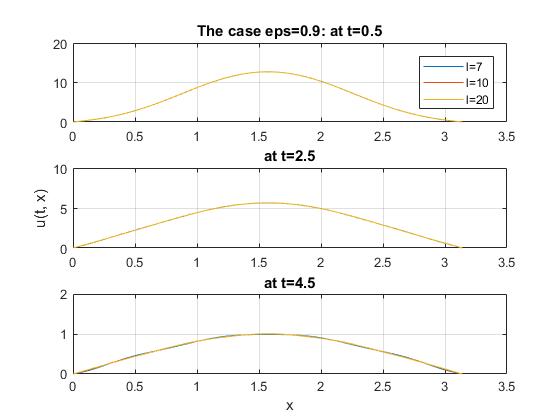}}
\end{minipage}
\caption{In the plots we illustrate the graphics of the solution $u$ of the inverse problem \eqref{CP-01}--\eqref{ED-01} for parameters $\varepsilon=0$ and $\varepsilon=0.9$ at $t=0.5, 2.5, 4.5$. Here the graphics of $f$ are given for different $l=7, 10, 20,$ and colored in blue, red, yellow, respectively.} \label{fig-1}
\end{figure}

For our calculations we use the following intermediate formulas:
\begin{equation}\label{Sol-num-test}
u(x,t)= \varphi(x) + \sum_{k=1}^{l} \frac{\int_{0}^{\pi} l^{2}(\psi(s) - \varphi(s)) v_k(s) ds}{\lambda_k^2 (1-e^{-\lambda_k T}) } (1 -  e^{-\lambda_k t}) v_k(x),
\end{equation}
and
\begin{equation}\label{f-num-test}
f(x)= l(\varphi(x)) + \sum_{k=1}^{l} \frac{\int_{0}^{\pi} l^{2}(\psi(s) - \varphi(s)) v_k(s) ds  }{ \lambda_k (1-e^{-\lambda_k T})} v_k(x).
\end{equation}
Here $l\in\mathbb N$ means a number of terms taken into account in the series of the formulas \eqref{Sol-num}--\eqref{f-num}.
Here, we analyse the convergence of the series in \eqref{Sol-num-test}--\eqref{f-num-test} for different $l$.

\begin{figure}[ht!]
\begin{minipage}[h]{0.49\linewidth}
\center{\includegraphics[scale=0.4]{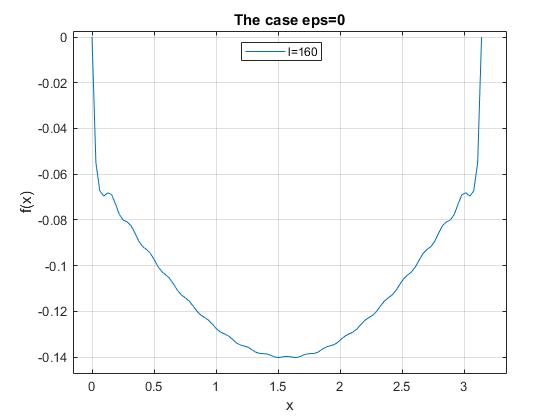}}
\end{minipage}
\hfill
\begin{minipage}[h]{0.49\linewidth}
\center{\includegraphics[scale=0.4]{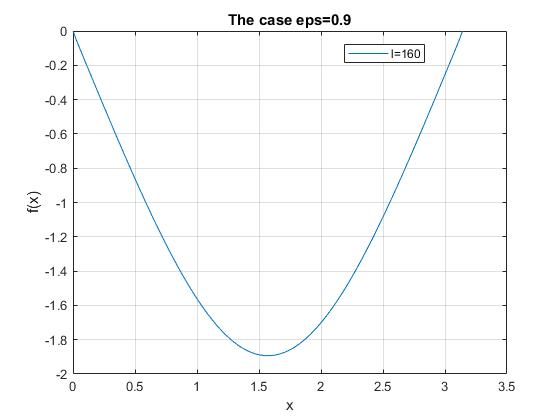}}
\end{minipage}
\caption{In these plots we compare the graphics of the source term $f$ of the inverse problem \eqref{CP-01}--\eqref{ED-01} for parameters $\varepsilon=0$ and $\varepsilon=0.9$.
} \label{fig3}
\end{figure}

\begin{figure}[ht!]
\begin{minipage}[h]{0.49\linewidth}
\center{\includegraphics[scale=0.4]{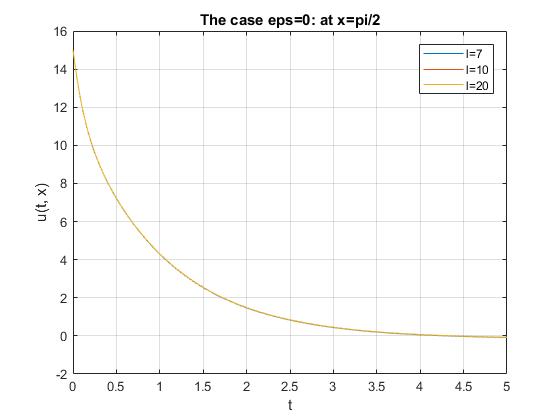}}
\end{minipage}
\hfill
\begin{minipage}[h]{0.49\linewidth}
\center{\includegraphics[scale=0.4]{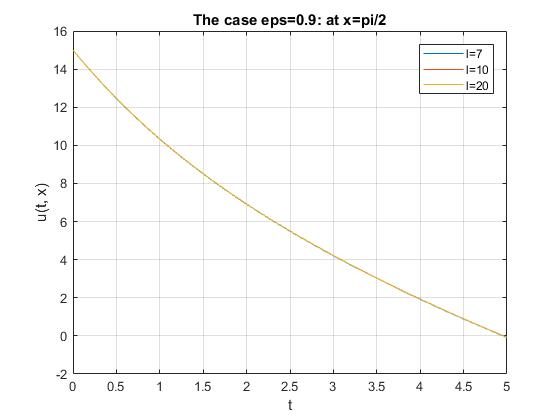}}
\end{minipage}
\caption{In the plots the graphics of the solution $u$ of the inverse problem \eqref{CP-01}--\eqref{ED-01} are shown at $x=\pi/2$ in the cases $\varepsilon=0$ and $\varepsilon=0.9$. Here the graphics of $f$ are given for different $l=7, 10, 20,$ and colored in blue, red, yellow, respectively.
} \label{fig4}
\end{figure}

In Figures \ref{fig-1}, \ref{fig3}, \ref{fig4} we compare the pair of solutions $(u, f)$ of the inverse problem \eqref{CP-01}--\eqref{ED-01} in the cases $\varepsilon=0$ and $\varepsilon=0.9$.

\subsection{Conclusion}

By analysing Figures \ref{fig-1}, \ref{fig3}, \ref{fig4}, we observe that in the case of the problem with an involution (when $\varepsilon\neq0$) we need more energy (see, Figure \ref{fig3}) to cool the rod with the initial temperature shown in Figure \ref{fig1}.
In Figure \ref{fig3} we can see the ``cooling speed" of the rod. As the result, the absence of the involution guarantees relatively rapid cooling.




\begin{thebibliography}{9}
\bibitem[ABGM15]{ABGM15} L. D. Abreu, P. Balazs, M. de Gosson and Z. Mouayn. Discrete coherent
states for higher Landau levels. {\em Ann. Physics}, 363: 337--353,
2015.

\bibitem[Bar51]{B51} N.~K.~Bari.
\newblock Biorthogonal systems and bases in Hilbert space.
\newblock {\it Moskov. Gos. Univ. Uchenye Zapiski Matematika}, 148(4):69--107, 1951.

\bibitem[BC17]{BC17} K. Beauchard, P. Cannarsa. Heat equation on the Heisenberg group: Observability and applications.
{\it J. Diff. Eq.} 262(8): 4475--4521, 2017.

\bibitem[CDR18]{CDR18}
M. Chatzakou, J. Delgado and M. Ruzhansky. On a class of anharmonic oscillators.
https://arxiv.org/abs/1811.12566

\bibitem[CS17]{CS17} L.~A.~Caffarelli, Y.~Sire. On Some Pointwise Inequalities Involving Nonlocal Operators. In: Chanillo S., Franchi B., Lu G., Perez C., Sawyer E. (eds) Harmonic Analysis, Partial Differential Equations and Applications. Applied and Numerical Harmonic Analysis. Birkh\"{a}user, Cham. 2017.

\bibitem[CD98]{CD98} J. R. Cannon, P. Du Chateau. Structural identification of an unknown source term in a heat equation. {\it Inverse Problems}. 14: 535--551, 1998.

\bibitem[CF18]{Carvalho}
P.~M.~de Carvalho-Neto, R. Fehlberg Junior.
Conditions for the Absence of Blowing Up Solutions to Fractional Differential Equations.
Acta Applicandae Mathematicae. 154(1): 15--29, 2018.

\bibitem[CNYY09]{CNYY09} J. Cheng, J. Nakagawa, M. Yamamoto, T. Yamazaki. Uniqueness in an inverse problem for a one-dimensional fractional diffusion equation. {\it Inverse Problems}. 25: 115002, 2009.

\bibitem[CG90]{CG90}
L.~J.~Corwin and F.~P.~Greenleaf.
\newblock {\em Representations of nilpotent Lie groups and their applications}.
\newblock Cambridge Studies in Advanced Mathematics, Cambridge University Press, Cambridge,
18. Basic theory and examples, 1990.

\bibitem[DRT17]{DRT17}
J. Delgado, M. Ruzhansky and N. Tokmagambetov.
\newblock Schatten classes, nuclearity and nonharmonic analysis on compact manifolds with boundary.
\newblock {\em J. Math. Pures Appl.}, 107:758--783, 2017.

\bibitem[DM05]{DM05} B. Driver, T. Melcher. Hypoelliptic heat kernel inequalities on the Heisenberg group. {\it J. Funct. Anal.} 221(5): 340--365, 2005.

\bibitem[FR16]{FR16}
V.~Fischer and M.~Ruzhansky.
\newblock {\em Quantization on nilpotent {L}ie groups}, volume 314 of {\em
  Progress in Mathematics}.
\newblock Birkh\"auser/Springer, [Open access book], 2016.

\bibitem[FR17]{FR:Sobolev}
V.~Fischer and M.~Ruzhansky.
\newblock Sobolev spaces on graded groups.
\newblock {\em Ann. Inst. Fourier}, 67(4):1671--1723, 2017.

\bibitem[Foc28]{F28} V. Fock, Bemerkung zur Quantelung des harmonischen
Oszillators im Magnetfeld. {\em Z. Phys. A}, 47(5--6):
446--448, 1928.

\bibitem[FS82]{FS-book}
G.~B. Folland and E.~M. Stein.
\newblock {\em Hardy spaces on homogeneous groups}, volume~28 of {\em
  Mathematical Notes}.
\newblock Princeton University Press, Princeton, N.J.; University of Tokyo
  Press, Tokyo, 1982.

\bibitem[FIK14]{Fur} K.M. Furati, O.S. Iyiola, M.~Kirane. An inverse problem for a generalized fractional diffusion. {\it Applied Mathematics and Computation}. 249: 24--31, 2014.

\bibitem[Gel63]{G63} I.~M.~Gelfand.
\newblock Some questions of analysis and differential equations.
\newblock {\it Am. Math. Soc. Transl.}, 26:201-219, 1963.

\bibitem[HH13]{HH13}
A.~Haimi and H.~Hedenmalm.
\newblock The polyanalytic Ginibre ensembles.
\newblock {\em J. Stat. Phys.}, 153(1):10--47, 2013.

\bibitem[HN79]{HN-79}
B.~Helffer and J.~Nourrigat.
\newblock Caracterisation des op\'erateurs hypoelliptiques homog\`enes
  invariants \`a gauche sur un groupe de {L}ie nilpotent gradu\'e.
\newblock {\em Comm. Partial Differential Equations}, 4(8):899--958, 1979.

\bibitem[HR82]{HR82}
B. Helffer and D. Robert.
Asymptotique des niveaux d'\'energie pour des hamiltoniens \`a un degr
\'e de libert\'e. {\em  Duke Math. J.}, 49(4):853--868, 1982.

\bibitem[HJL85]{HJL85}
A.~Hulanicki, J.~W.~Jenkins, and J.~Ludwig.
\newblock Minimum eigenvalues for positive, Rockland operators.
\newblock {\em Proc. Amer. Math. Soc.}, 94:718--720, 1985.

\bibitem[IC16]{Ismailov} M. I. Ismailov, M. Cicek. Inverse source problem for a time-fractional diffusion equation with nonlocal boundary conditions. {\it Applied Mathematical Modelling}. 40(7): 4891--4899, 2016.

\bibitem[JR15]{JR15} B. Jin, W. Rundell. A tutorial on inverse problems for anomalous diffusion processes. {\it Inverse Problems}. 31(3): 035003, 2015.

\bibitem[Jui14]{Jui14} N. Juillet. Diffusion by optimal transport in Heisenberg groups. {\it Calc. Var. Partial Differ. Equ.} 50(3-4): 693--721, 2014.

\bibitem[KS10]{Kali1} I. A. Kaliev, M.M. Sabitova. Problems of determining the
temperature and density of heat sources from the initial and final
temperatures. {\it Journal of Applied and Industrial Mathematics}. 4(3): 332--339, 2010.

\bibitem[KRT17]{KRT17}
B.~Kanguzhin, M.~Ruzhansky, and N.~Tokmagambetov.
\newblock On convolutions in Hilbert spaces.
\newblock {\em Funct. Anal. Appl.}, 51(3):221--224, 2017.

\bibitem[KST06]{Kilbas} A. A. Kilbas, H. M. Srivastava, J.J. Trujillo. \emph{Theory and Applications of Fractional Differential Equations}, Elsevier, North-Holland, Mathematics studies, 2006.

\bibitem[KM11]{Kir} M. Kirane, A.S. Malik. Determination of an unknown source term and the temperature distribution for the linear heat equation involving fractional derivative in time. {\it Applied Mathematics and Computation}. 218(1):163--170, 2011.

\bibitem[KST17]{Kir1} M. Kirane, B. Samet, B. T. Torebek. Determination of an unknown source term temperature distribution for the sub-diffusion equation at the initial and final data. {\it Electronic Journal of Differential Equations}. 2017: 1--13, 2017.

\bibitem[Lan30]{L30}
L. Landau.
Diamagnetismus der Metalle.
{\em Z. Phys. A}, 64(9--10): 629--637, 1930.

\bibitem[LG99]{LG99}
Y.~Luchko, R.~Gorenflo.
\newblock An operational method for solving fractional differential equations
with the Caputo derivatives.
\newblock {\it Acta Math. Vietnam.}, 24:207--233, 1999.

\bibitem[Mai00]{Mainardi} F. Mainardi. Waves and Stability in Continuous Media. Singapore: World Scientific, 2000.

\bibitem[MRT19a]{MRT19a}
J.~C.~Munoz, M.~Ruzhansky, and N.~Tokmagambetov.
\newblock Acoustic and Shallow Water Wave Propagation with Irregular Dissipation.
\newblock {\em Functional analysis and its applications}, 53(2): 153-156, 2019.

\bibitem[MRT19b]{MRT19b}
J.~C.~Munoz, M.~Ruzhansky, and N.~Tokmagambetov.
\newblock Wave propagation with irregular dissipation and applications to acoustic problems and shallow waters.
\newblock {\em Journal de mathematiques pures et appliquees}, 123: 127-147, 2019.

\bibitem[Nai68]{Naimark}
M.~A. Na{\u\i}mark.
\newblock {\em Linear differential operators. {P}art {II}: {L}inear
  differential operators in {H}ilbert space}.
\newblock With additional material by the author, and a supplement by V. {\`E}.
  Ljance. Translated from the Russian by E. R. Dawson. English translation
  edited by W. N. Everitt. Frederick Ungar Publishing Co., New York, 1968.

\bibitem[NLN16]{Nguyen} H. T. Nguyen, D. L. Le, V. T. Nguyen. Regularized solution of an inverse source problem for a time fractional diffusion equation. {\it Applied Mathematical Modelling}. 40(19): 8244--8264, 2016.

\bibitem[NR10]{NiRo:10}
F.~Nicola and L.~Rodino.
\newblock {\em Global pseudo-differential calculus on {E}uclidean spaces},
  volume~4 of {\em Pseudo-Differential Operators. Theory and Applications}.
\newblock Birkh{\"a}user Verlag, Basel, 2010.

\bibitem[Nig86]{Nigmatullin} R. R. Nigmatullin. The realization of the generalized transfer equation in a medium with
fractal geometry, {\it Phis. Stat. Sol.}, 133:299--318, 1986.

\bibitem[OS12a]{Oraz} I. Orazov, M.A. Sadybekov. One nonlocal problem of
determination of the temperature and density of heat sources.
{\it Russian Mathematics}. 56(2):60--64, 2012.

\bibitem[OS12b]{Oraz1} I. Orazov, M.A. Sadybekov. On a class of problems of
determining the temperature and density of heat sources given
initial and final temperature. {\it Siberian Mathematical Journal}. 53(1):146--151, 2012.

\bibitem[RS80]{RS80}
M.~Reed and B.~Simon.
\newblock {\em Methods of Modern Mathematical Physics}, volume~1 of {\em
  Functional Analysis,
revised and enlarged edition}.
\newblock Academic Press, 1980.

\bibitem[Roc78]{Rockland}
C.~Rockland.
\newblock Hypoellipticity on the {H}eisenberg group-representation-theoretic
  criteria.
\newblock {\em Trans. Amer. Math. Soc.}, 240:1--52, 1978.

\bibitem[RR18]{RR18}
D. Rottensteiner, M. Ruzhansky.
\emph{Harmonic and Anharmonic Oscillators on the Heisenberg group}. https://arxiv.org/abs/1812.09620

\bibitem[RS76]{RS76}
L.~P.~Rothschild and E.~M.~Stein.
\newblock Hypoelliptic differential operators and nilpotent
groups.
\newblock {\em Acta Math.}, 137:247--320, 1976.

\bibitem[RT16]{RT16}
M.~Ruzhansky, N.~Tokmagambetov.
\newblock Nonharmonic analysis of boundary value problems.
\newblock {\it Int. Math. Res. Notices}, 2016(12): 3548--3615, 2016.

\bibitem[RT17a]{RT17a}
M.~Ruzhansky, N.~Tokmagambetov.
\newblock Nonharmonic analysis of boundary value problems without WZ condition.
\newblock {\it Math. Model. Nat. Phenom.}, 12(1):115--140, 2017.

\bibitem[RT17b]{RT17b}
M.~Ruzhansky, N.~Tokmagambetov.
\newblock Very weak solutions of wave equation for Landau Hamiltonian with irregular electromagnetic field.
\newblock {\it Lett. Math. Phys.}, 107:591--618, 2017.

\bibitem[RT18]{RT18}
M.~Ruzhansky and N.~Tokmagambetov.
\newblock On a very weak solution of the wave equation for a Hamiltonian in a singular electromagnetic field.
\newblock {\em Math. Notes}, 103(5--6):856--858, 2018.

\bibitem[RT18b]{RT18b}
M.~Ruzhansky and N.~Tokmagambetov.
\newblock Convolution, Fourier analysis, and distributions generated by Riesz bases.
\newblock {\em Monatsh. Math.}, 187(1):147–-170, 2018.

\bibitem[RT17c]{RT17c}
M.~Ruzhansky, N.~Tokmagambetov.
\newblock Wave equation for operators with discrete spectrum and irregular propagation speed.
\newblock {\it Arch. Ration. Mech. Anal.}, 226:1161--1207, 2017.

\bibitem[RT18c]{RT18c}
M.~Ruzhansky and N.~Tokmagambetov.
\newblock Nonlinear damped wave equations for the sub--Laplacian on the Heisenberg group and for Rockland operators on graded Lie groups.
\newblock {\em Journal of Differential Equations}, 265(10):5212--5236, 2018.

\bibitem[RTT19]{RTT19}
M.~Ruzhansky, N.~Tokmagambetov, and B.T.~Torebek.
\newblock Bitsadze-Samarskii type problem for the integro-differential diffusion-wave equation on the Heisenberg group,
\newblock {\em Integral transforms and special functions}, 2019, to appear.

\bibitem[RT19a]{RT19a}
M.~Ruzhansky and N.~Tokmagambetov.
\newblock Wave Equation for 2D Landau Hamiltonian.
\newblock {\em Applied and computational mathematics}, 18(1): 69-78, 2019.

\bibitem[RT19b]{RT19b}
M.~Ruzhansky and N.~Tokmagambetov.
\newblock On nonlinear damped wave equations for positive operators. I. Discrete spectrum. \newblock {\em Differential and integral equations}, 32(7-8): 455-478, 2019.

\bibitem[SY11]{SY11} K. Sakamoto, M. Yamamoto. Inverse source problem with a final overdetermination for afractional diffusion equation. {\it Math. Control Relat. Fields}. 1: 509--518, 2011.

\bibitem[Sim14]{Simon}
T.~Simon.
Comparing Frechet and positive stable laws. {\it Electron. J. Probab}. 19:1--25, 2014.

\bibitem[TT16]{TT16} N. Tokmagambetov, B.T. Torebek.
Fractional Analogue of Sturm--Liouville Operator. {\it Documenta Mathematica}. 21: 1503--1514, 2016.

\bibitem[TT18]{TT18} N. Tokmagambetov, B.T. Torebek.
Well--posed problems for the fractional Laplace equation with integral boundary conditions.
{\it Electronic Journal of Differential Equations}. 2018: 1--10, 2018.

\bibitem[TT18d]{TT18d}
N.~Tokmagambetov, B.~T.~Torebek.
\newblock Green's formula for integro--differential operators.
\newblock {\em J. Math. Anal. Appl.}, 468(1):473--479, 2018.

\bibitem[TT18a]{TT18a}
Tokmagambetov N., Torebek B.T.
Fractional Sturm--Liouville Equations: Self--Adjoint Extensions.
{\it Complex Analysis and Operator Theory}, 13(5):2259--2267, 2019.

\bibitem[TT17]{torebek} B.T. Torebek, R. Tapdigoglu. Some inverse problems for the nonlocal heat equation with Caputo fractional derivative. {\it Mathematical Methods in the Applied Sciences}. 40(18): 6468--6479, 2017.

\bibitem[TW18]{TW18} J. Tyson, J. Wang. Heat content and horizontal mean curvature on the Heisenberg group. {\it Communications in Partial Differential Equations}. 43(3): 467--505, 2018.

\bibitem[Uch13]{Uchaikin}
V.V. Uchaikin. \emph{Fractional Derivatives for Physicists and Engineers}. V. 1, Background and Theory. V. 2, Application.
Springer, 2013.

\bibitem[ZX11]{ZX11} Y. Zhang, X. Xu. Inverse source problem for a fractional diffusion equation. {\it Inverse Problems}. 27: 035010, 2011.

\bibitem[ZW13]{ZW13} Z. Q. Zhang, T. Wei. Identifying an unknown source in time-fractional diffusion equation by a truncation method. {\it Appl. Math. Comput.} 219: 5972--5983, 2013.

\bibitem[WYH13]{WYH13} W. Wang, M. Yamamoto, B. Han. Numerical method in reproducing kernel space for an inverse source problem for the fractional diffusion equation. {\it Inverse Problems}. 29: 095009, 2013.

\bibitem[Wey31]{W31} H. Weyl, The Theory of Groups and Quantum Mechanics, Methuen, 1931.

\end{thebibliography}
\end{document}